\documentclass{siamart}

\usepackage{lipsum}
\usepackage{amsfonts}
\usepackage{graphicx}
\usepackage{amsopn}

\usepackage{enumerate}

\usepackage{mathtools}
\usepackage{booktabs}

\usepackage{color}
\usepackage{changes}

\usepackage{epstopdf}
\usepackage{algorithmic}

\newsiamremark{remark}{Remark}
\newtheorem*{example*}{Example}
\newtheorem{conjecture}{Conjecture}

\newcommand{\N}{\mathbb N}

\newcommand{\R}{\mathbb R}
\newcommand{\C}{\mathbb C}

\newcommand{\mcM}{\mathcal M}

\newcommand{\mcG}{\mathcal G}
\newcommand{\mcH}{\mathcal H}

\newcommand{\colspan}{\operatorname{ran}}
\newcommand{\rank}{\operatorname{rank}}
\newcommand{\tr}{\operatorname{tr}}
\newcommand{\sym}{\operatorname{sym}}
\newcommand{\Skew}{\operatorname{skew}} 

\newcommand{\Ad}{\operatorname{Ad}}

\newcommand{\diff}[1][t]{\mathrm{d}}
\newcommand{\Diff}[1][t]{\mathrm{D}}

\makeatletter
\renewcommand*\env@matrix[1][*\c@MaxMatrixCols c]{%
	\hskip -\arraycolsep
	\let\@ifnextchar\new@ifnextchar
	\array{#1}}
\makeatother

\newcount\Comments  
\newcommand{\kibitz}[2]{\ifnum\Comments=1\textcolor{#1}{#2}\fi}

\newcommand{\RZ}[1]  {\kibitz{teal}   {[RZ: #1]}}
\newcommand{\JS}[1]  {\kibitz{blue}   {[JS: #1]}}

\ifpdf
\DeclareGraphicsExtensions{.eps,.pdf,.png,.jpg}
\else
\DeclareGraphicsExtensions{.eps}
\fi

\newcommand{\TheTitle}{
	High curvature means low-rank: On the sectional curvature of Grassmann and Stiefel manifolds and the underlying matrix trace inequalities}
\newcommand{\TheAuthors}{Ralf Zimmermann and Jakob Stoye}

\headers{High curvature means low-rank}{\TheAuthors}

\title{{\TheTitle}} 

\author{
Ralf Zimmermann%
\thanks{University of Southern Denmark, Department of Mathematics and Computer Science, Odense, Denmark
	(\email{zimmermann@imadasdu.dk}, \url{https://portal.findresearcher.sdu.dk/en/persons/zimmermann}, 
	\url{https://orcid.org/0000-0003-1692-3996})
}
\and
Jakob Stoye%
\thanks{TU Braunschweig, Institute for Numerical Analysis, Braunschweig, Germany,
	(\email{jakob.stoye@tu-braunschweig.de}, \url{https://orcid.org/0009-0003-9119-0013})
}
}

\ifpdf
\hypersetup{
  pdftitle={\TheTitle},
  pdfauthor={\TheAuthors}
}
\fi

\begin{document}

\maketitle

\begin{abstract}
Methods and algorithms that work with data on nonlinear manifolds are collectively summarized under the term `Riemannian computing'.
In practice, curvature can be a key limiting factor for the performance of Riemannian computing methods. Yet, curvature can also be a powerful tool in the theoretical analysis of Riemannian algorithms.
In this work, we investigate the sectional curvature of the Stiefel and Grassmann manifold.
On the Grassmannian, tight curvature bounds are known since the late 1960ies. On the Stiefel manifold under the canonical metric, it was believed that the sectional curvature does not exceed 5/4. Under the Euclidean metric, the maximum was conjectured to be at 1.
For both manifolds, the sectional curvature is given by the Frobenius norm of certain structured commutator brackets of skew-symmetric matrices. We provide refined inequalities for such terms and pay special attention to the maximizers of the curvature bounds. In this way, we prove for the Stiefel manifold that the global bounds of 5/4 (canonical metric) and 1 (Euclidean metric) hold indeed. With this addition, a complete account of the curvature bounds in all admissible dimensions is obtained.
We observe that `high curvature means low-rank', more precisely, for the Stiefel and Grassmann manifolds under the canonical metric, the global curvature maximum is attained at tangent plane sections that are spanned by rank-two matrices, while the extreme curvature cases of the Euclidean Stiefel manifold occur for rank-one matrices.
Numerical examples are included for illustration purposes.
\end{abstract} 


\section{Introduction}
\label{sec:intro}
Methods and algorithms that work with data on nonlinear manifolds are collectively summarised under the term Riemannian computing.
Riemannian computing methods have established themselves as important tools in a large variety of applications, including  computer vision, machine learning, and optimization, see
\cite{AbsilMahonySepulchre2004,AbsilMahonySepulchre2008,BegelforWerman2006,EdelmanAriasSmith1999,Gallivan_etal2003,  Lui2012, Rahman_etal2005, Rentmeesters2013} and the anthologies \cite{Minh:2016:AAR:3029338, RiemannInComputerVision}.
They also gain increasing attention in statistics and data science
\cite{patrangenaru2015nonparametric} 
and in numerical methods for differential equations \cite{BennerGugercinWillcox2015,Celledoni_2020,hairer06gni, iserles_munthe-kaas_norsett_zanna_2000,Zimmermann_MORHB2021}.\\
A standard technique in designing Riemannian computing methods is to translate Euclidean algorithms to manifolds.
The discrepancy between a linear space and a non-linear manifold is quantified by the concept of curvature. Therefore, curvature can also be seen as a decisive factor that separates Riemannian algorithms from their Euclidean counterparts.
On the other hand, curvature estimates allow to compare Euclidean distances and Riemannian distances for embedded submanifolds \cite{AttaliEM07}, which can be a powerful tool in the theoretical analysis of Riemannian methods.

In this work, we investigate the sectional curvature on the Stiefel and Grassmann manifolds.
Both of these manifolds have been subject to extensive investigations before.
On the Grassmannian, tight curvature bounds are known since the seminal papers of Wong \cite{Wong1967, Wong1968}. 
The sectional curvature on the Stiefel manifold under a parametric family of Riemannian metrics introduced in \cite{HueperMarkinaLeite2020} was extensively studied in \cite{nguyen2022curvature}.
On the Stiefel manifold under the canonical metric,
it was believed that the sectional curvature does not exceed $5/4$, \cite[p. 94,95]{Rentmeesters2013}, \cite[Section 6, Table 2]{nguyen2022curvature}. 
Propositions. 6.1 and 6.2 from the same work report that under the Euclidean metric, the curvature range contains the interval $[-1/2, 1]$.
Nguyen \cite{nguyen2022curvature} shows that these bounds are tight for the Stiefel manifold $St(n,p)$ if  $p=2$, but a valid proof for general $p\geq2$ was missing.

For the orthogonal group and, in turn, for the Stiefel and Grassmann manifolds as its quotient spaces, the sectional curvature is given by the Frobenius norm of certain structured commutator brackets of skew-symmetric matrices. 
Sharp estimates are provided by the matrix inequality of \cite{WuChen1988} and the special bounds for skew-symmetric matrices of \cite[Lemma 2.5]{Ge2014}.
We provide refined inequalities for such terms and pay special attention to the maximizers of the curvature bounds. In this way, we prove the conjectured global curvature bounds of
\[
    0\leq K_c \leq \frac54  \text{ (canonical metric)} \quad \text{and} \quad 
    -\frac12\leq K_e \leq 1  \text{ (Euclidean metric)}.
\]
We also show that these bounds are sharp on all true Stiefel manifolds $St(n,p)$, where $p\geq 2, n\geq p+2$ and provide a discussion of the special cases, where $p=1$, $p=n$, and $p=n-1$.

In doing so, we observe that `high curvature means low-rank'. 
To be precise, we prove that for the orthogonal group and the Stiefel and Grassmann manifolds under the canonical metric, the global curvature maximum is attained at tangent plane sections that are spanned by 
the same low-rank tangent matrices.\footnote{Because Stiefel and Grassmann are quotient spaces of the orthogonal group $SO(n)$, in this case, Stiefel and Grassmann tangent vectors can be considered as special (namely, horizontal) $SO(n)$-tangent vectors.}
From the perspective of the orthogonal group, these are rank-four matrices, while they are
rank-two matrices from the Stiefel or Grassmann view point. 
Under the Euclidean metric, the extreme curvature cases occur for rank-one matrices.
Numerical experiments confirm that the curvature drops with increasing rank.

\paragraph{Organisation of the paper}
\Cref{sec:LieGroups} provides the required curvature formulas for matrix Lie groups and their quotient spaces.
(For the reader's convenience, some essentials of Lie group theory
is gathered in \Cref{app:matLieGroups}.)
Readers not interested in the matrix manifold applications but only in the matrix norm inequalities may directly skip to \Cref{sec:curv_estimates}.
In this section, we also give a full account of the sharp Stiefel curvature bounds in all possible dimensions. Numerical experiments are discussed in \Cref{sec:numex}.
\Cref{sec:conclusions} concludes the paper.
\paragraph{Notational specifics}
\label{sec:Notation}
For $n\in \N$, the $(n\times n)$-identity matrix is denoted by $I_n\in\R^{n\times n}$, or simply $I$, if the dimension is clear.
The $(n\times n)$ special orthogonal group, i.e., the set of all
square orthogonal matrices with determinant $1$ is denoted by
\[
  SO(n) = \{Q\in \R^{n\times n}| Q^TQ = QQ^T = I_n\}.
\]
The orthogonal group $O(n)$ is $SO(n)$ joint with their `$\det=-1$'-siblings.
The sets of symmetric and skew-symmetric $(n\times n)$-matrices are $\sym(n) = \{A\in\R^{n\times n}|A^T=A\}$ and
$\Skew(n) = \{A\in\R^{n\times n}|A^T=-A\}$, respectively.
Overloading this notation, $\sym(A) = \frac12(A+A^T), \Skew(A) = \frac12(A-A^T)$
denote the symmetric and skew-symmetric parts of a matrix $A$.
The Stiefel manifold is
\[
    St(n,p) = \{U\in\R^{n\times p}\mid U^TU=I_p\}.
\]
The Grassmann manifold is
\[
    Gr(n,p) = \{[U]\mid U\in St(n,p)\}, \text{where } [U] = \{UR\mid R\in SO(p)\}.
\]
{\em A word of caution:}
The standard Euclidean inner product on $\R^{n\times p}$ is 
$$\langle X,Y\rangle_F = \tr(X^TY).$$
The subscript is to emphazise the correspondence with the Frobenius norm $\|X\|_F = \sqrt{\langle X,X\rangle_F}$. To comply with standard conventions in the 
matrix manifolds and Riemannian computing literature, the Riemannian metric on $SO(n)$ will \emph{not}
be the Euclidean one, but the Euclidean one with a multiplicative factor of $\frac12$, i.e.,
\begin{equation}
	\label{eq:standard_metric}
	\langle X,Y\rangle_Q = \frac12 \tr(X^TY)= \frac12 \langle X,Y\rangle_F, \quad X,Y \in T_QSO(n).
\end{equation}
This factor is inherited by the Riemannian metrics of the Grassmann and the Stiefel manifold,
when considered as quotients of $SO(n)$.
It also makes the curvature results in this work compatible with those stated, e.g., in 
the seminal papers \cite{Wong1967,Wong1968}.
To distinguish the Riemannian and the Euclidean metric, for the former, the location is always given as a subscript $\langle\cdot,\cdot\rangle_Q$.

Let $\mcM$ be a quotient of $SO(n)$ under a Lie group action. (We will only consider $\mcM=St(n,p)$ or $\mcM=Gr(n,p)$.)
For a tangent vector $X\in T_ISO(n)= \mathfrak{so}(n)=\Skew(n)$, 
$X_{\mathfrak{m}}$ denotes the projection onto the horizontal space associated with $\mcM$.
	When a distinction is necessary, we will also write $X_{\mathfrak{m}^\mcM}$ to emphasize, which quotient manifold is considered.
Likewise, $X_{\mathfrak{h}}$ and $X_{\mathfrak{h}^\mcM}$ denote the projection onto the vertical space associated with $\mcM$.
%
%
%
\section{Curvature of Lie groups and quotients of Lie groups}
\label{sec:LieGroups}
This section recaps the basic curvature formulas for Lie groups and quotients of Lie groups. This is classic textbook material. Our main references are  \cite{gallier2020differential}, \cite{Hall_Lie2015} and  \cite{ONeill1983}, to which we refer for the details. Readers not familiar with Lie groups or quotients of Lie groups may want to read \Cref{app:LieGroups} first.
%
\subsection{Lie group curvature formulae}
Let $\mathcal{G}$ be a Lie group with Lie algebra $\mathfrak{g}=T_I\mathcal{G}$ (the tangent space at $I$) and a bi-invariant metric.
Let $X,Y\in \mathfrak{g}$ be linearly independent tangent vectors. The sectional curvature associated with the two-dimensional subplane spanned by $\{X,Y\}\subset \mathfrak{g}$ is
\[
    K(X,Y) = \frac14\frac{\|[X,Y]\|_I^2}{\|X\|_I^2\|Y\|_I^2-\langle X,Y\rangle_I},
\]
where $[X,Y] = XY-YX$.
It depends only on the subplane, in this context often referred to as `{\em the section}', not on the spanning vectors $X,Y$.
For convenience, we will only consider $X,Y\in \mathfrak g$ forming an orthonormal basis (ONB),
i.e., $\|X\|_I = 1 = \|Y\|_I$, $\langle X,Y\rangle_I =0$, so that the
sectional curvature is computed as
\begin{equation}
\label{eq:curv_biinv_Liegroup}
	K(X,Y) = \frac14 \|[X,Y]\|_I^2,
\end{equation}
see \cite[Prop. 21.19, p. 636]{gallier2020differential}.\\
\begin{remark}
\label{rem:curvSOn_fro}
For $\mathcal{G}=SO(n)$, we use the metric of \eqref{eq:standard_metric}.
Hence, when translated to the Frobenius norm, one has
\[
	K^{SO}(X,Y) = \frac12 \|[X,Y]\|_F^2,\quad \text{for } \|X\|_F=1=\|Y\|_F, \quad \langle X,Y\rangle_F =0.
\]
\end{remark}
In order to obtain curvature expressions for Lie group quotients,
we need to introduce a few more terms.
The {\em adjoint map} ${\bf Ad}: \mcG \to GL(\mcG)$ assigns to each group element $A\in\mcG$ the linear map ${\bf Ad}_A: B\mapsto ABA^{-1}$. Note that ${\bf Ad}_A(I) = I$.
The differential of ${\bf Ad}_A$ at $I$ is the map $\Ad_A:=d({\bf Ad}_A)_I:T_I\mcG = \mathfrak g \to  \mathfrak g=T_{{\bf Ad}_A(I)}\mcG$. For a matrix Lie group, 
$\Ad_A(X) = AXA^{-1}$.

The homogeneous $\mcG$-space $\mcG/\mcH$ is {\em reductive},
if the Lie algebra $\mathfrak g$ (the tangent space at the identity)
can be split into 
\[
		\mathfrak g = \mathfrak h \oplus \mathfrak m,
\]
where $\mathfrak h$ is the Lie algebra of $\mcH$ (the tangent space at the identity) and $\mathfrak m\subset \mathfrak g$ is a complementary subspace
such that $\Ad_Q(\mathfrak m)\subseteq \mathfrak m$ for all $Q\in \mcH$.
Note that while the setting above is more general, in the cases that we 
eventually consider in this work, $\mathfrak m$ will always be the {\em orthogonal} complement of $\mathfrak h$ with respect to a suitable metric.
For a tangent vector $X\in\mathfrak g$, the projections onto
 $\mathfrak h$ and $\mathfrak m$ are denoted by $X_{\mathfrak h}$
 and $X_{\mathfrak m}$, respectively.

The quotient space $\mcG/\mcH$ is called {\em naturally reductive}, if
\begin{enumerate}[(i)]
	\item it is reductive with some decomposition $\mathfrak g = \mathfrak h \oplus \mathfrak m$,
	\item it has a $\mcG$-invariant metric,
	\item $\langle [X,Z]_{\mathfrak m},Y\rangle_I = \langle X,[Z,Y]_{\mathfrak m}\rangle_I$ for all $X,Y,Z\in \mathfrak m$.
\end{enumerate}
Item (ii) means that the maps $\tau_G: \mcG/\mcH \to \mcG/\mcH,[A]\mapsto G\cdot[A]$ associated with the action $\mcG\times \mcG/\mcH \to \mcH, (G, [A])\mapsto G\cdot [A]$ are isometries.
The precise condition is
\[
\langle d(\tau_G)_{[A]} (X), d(\tau_G)_{[A]} (Y)\rangle_{\tau_G([A])}=\langle X,Y\rangle_{[A]} \quad \forall X,Y\in T_{[A]}\mcG/\mcH, [A]\in \mcG/\mcH,
\]
 \cite[Def. 23.5]{gallier2020differential}.
On the level of matrix representatives (and after identifying $T_{[A]}\mcG/\mcH$ with the horizontal space $H_A$),
this boils down to the condition
\[
   \langle G\bar X,G\bar Y\rangle_{GA} =\langle \bar X,\bar Y\rangle_A \text{ for all } \bar X, \bar Y\in  H_A, A\in \mcG,
\]
where $\bar X, \bar Y\in H_A$ are horizontal lifts of $X,Y\in T_{[A]}\mcG/\mcH$.
	\begin{theorem}[\cite{gallier2020differential}, Prop. 23.29]
 \label{thm:quotcurv}
	Let  $\mcG$ be a connected 	Lie group such that the Lie algebra $\mathfrak g$ 
	admits an $\Ad(\mcG)$-invariant inner product $\langle \cdot,\cdot\rangle_I$. 
	Let $\mcG/\mcH$ be a homogeneous space as above 
	and let $\mathfrak m =  \mathfrak{h}^\bot$  with respect to $\langle \cdot,\cdot\rangle_I$. 
	Then
	\begin{enumerate}
			\item The space $\mcG/\mcH$ is reductive with respect to the decomposition $\mathfrak g =\mathfrak h \oplus \mathfrak m$.
			\item Under the $\mcG$-invariant metric induced by the inner product, the homogeneous space $\mcG/\mcH$
			is naturally reductive.
		\item For an ONB spanned by $X,Y\in T_{[I]}\mcG/\mcH\cong \mathfrak m$, $X\bot Y$, $\|X\|_I=\|Y\|_I=1$, the sectional curvature 
		associated with the tangent plane spanned by $X,Y$ is
		\begin{equation}
			\label{eq:curv_basis}
		 K(X,Y) = 
			\frac14 \| [X,Y]_{\mathfrak m} \|^2_I + \| [X,Y]_{\mathfrak h} \|_I^2.
		\end{equation}
	\end{enumerate}
\end{theorem}
\paragraph{Transport to arbitrary locations}
Since the tangent space at an arbitrary location $A\in\mcG$ is given by the translates of the tangent space at the identity (see \eqref{eq:tangspaceshifts}), a tangent vector $\tilde X\in T_A\mcG$ at an arbitrary location $A\in \mcG$ is of the form $\tilde X = AX$ with $X\in \mathfrak g$.
For an ONB $\{\tilde X=AX, \tilde Y=AY\}$ that spans a tangent plane in $T_{A}\mcG$,
the sectional curvature is
\[
K^{\mathcal{G}}(\tilde X,\tilde Y) = \frac14 \|[X,Y]\|_I^2.
\]

For quotient spaces, because of the natural reductive homogeneous space structure,
the horizontal spaces $H_A$ at arbitrary locations $A\in \mcG$ are translates of the horizontal space 
$\mathfrak m \cong T_{[I]}\mcG/\mcH$. Therefore, the tangent space of the quotient $\mcG/\mcH$ is identified with
\[
	T_{[A]}\mcG/\mcH \cong H_A = A\mathfrak m.
\]
Let $\tilde X,\tilde Y\in T_{[A]}\mcG/\mcH$ be an ONB of a tangent plane in the quotient space.
Let $\bar X, \bar Y\in H_A\subset T_A\mcG$ be horizontal lifts of 
$\tilde X,\tilde Y\in T_{[A]}\mcG/\mcH \cong A\mathfrak m$ with $\bar X = AX, \bar Y = AY$, $X,Y\in \mathfrak{m}$.
Then the sectional curvature 
at $[A]\in T_{[A]}\mcG/\mcH$ with respect to the tangent plane is given by
\begin{equation}
	K^{\mcG/\mcH}(\tilde X,\tilde Y) = 
	\frac14 \| [X,Y]_{\mathfrak m} \|^2_I + \| [ X, Y]_{\mathfrak h} \|^2_I.
\end{equation}
For details, see \cite[\S 19--23]{gallier2020differential}.

%
%
%
\subsection{Canonical curvature formulae on Grassmann and Stiefel}
In this section we recap the expressions for the sectional curvature on the Grassmann and the Stiefel manifold.
Both the Stiefel manifold $St(n,p)$ and the Grassmann manifold $Gr(n,p)$ are considered as quotients of $SO(n)$.
The Riemannian metric on the total space $SO(n)$ is 
\begin{equation}
	\label{eq:innerSOn}
	\langle X, Y\rangle_Q = \frac12 \tr(X^T Y), \quad X,Y\in T_QSO(n) = Q\Skew(n).
\end{equation}
All conditions of \Cref{thm:quotcurv} are fulfilled for $\mathcal{G}=SO(n)$ with this metric and its quotient spaces $St(n,p)$ and $Gr(n,p)$.
\paragraph{Grassmann sectional curvature}
The Grassmann manifold of linear subspaces can be realized as 
the quotient space $Gr(n, p) = SO(n)/(SO(p)\times SO(n-p))$.
The canonical projection and the left cosets are 
\begin{eqnarray*}
	\label{eq:Grassmann_quotient1}
	\Pi^{Gr}&:&SO(n)\to Gr(n,p)\\ 
	\label{eq:Grassmann_quotient_cosets} 
	[Q]=: \Pi^{Gr}(Q)  &=& \{Q\begin{pmatrix}
		S& 0\\
		0 & R
	\end{pmatrix}\mid \begin{pmatrix}
	S& 0\\
	0 & R
\end{pmatrix}\in SO(p)\times SO(n-p)\}.
\end{eqnarray*}
By splitting $Q = \begin{pmatrix}[c|c]
	U & U_\bot
\end{pmatrix}$ with $U\in\R^{n\times p}, U_\bot\in \R^{n\times (n-p)}$,    $[Q]\in Gr(n,p)$ is uniquely represented by $\colspan(U)\subset \R^{n}$.\\
When lifting $[Q]\in Gr(n,p)$ to $Q\in SO(n)$ (in practice, this is nothing but fixing a representative $Q$ for the equivalence class $[Q]$), the vertical space is represented by
		\begin{equation*}
			V_Q^{Gr} = 
			\bigl\{ \bar X = Q
			\begin{pmatrix}
				A & 0\\
				0 & C
			\end{pmatrix}\mid 
			A\in \Skew(p), C\in \Skew(n-p)
			\bigr\}=Q\mathfrak{h}^{Gr}.
		\end{equation*}
The associated horizontal space is
		\begin{equation*}
		T_{[Q]}Gr(n,p) \cong 	H_Q^{Gr} =
			\bigl\{\bar X = Q
			\begin{pmatrix}
				0 & -B^T\\
				B &  0
			\end{pmatrix}\mid 
			 B\in \R^{(n-p)\times p}
			\bigr\} = Q\mathfrak{m}^{Gr}.
		\end{equation*}
For Grassmann tangent vectors (already represented in matrix form by their lifts) $ X =\begin{pmatrix}
	0 & -B^T_1\\
	B_1 &  0
\end{pmatrix},
  Y=			\begin{pmatrix}
 	0 & -B^T_2\\
 	B_2 &  0
 \end{pmatrix}\in \mathfrak m^{Gr} \cong T_{[I]}Gr(n,m)$,
it holds that $ X, Y$ are orthonormal w.r.t. \eqref{eq:innerSOn} if and only if
$B_1,B_2\in \R^{(n-p)\times p}$ are orthonormal w.r.t. the standard Euclidean inner product
$\langle \cdot,\cdot\rangle_F$.
For such orthonormal $X,Y$, a straightforward evaluation of \eqref{eq:curv_basis} yields
\begin{eqnarray}
		\nonumber
		K^{Gr}(X,Y) &=&	\frac14 \| [ X, Y]_{\mathfrak{m}^{Gr}} \|^2_{[I]} + \| [ X, Y]_{\mathfrak{h}^{Gr}} \|^2_{[I]}\\	
		\label{eq:seccurv_Grassmann1}
		&=& \frac12 \|B_1^TB_2 - B_2^TB_1\|_F^2 +  \frac12\|B_1B_2^T - B_2B_1^T\|_F^2\\  
	  \label{eq:seccurv_Grassmann2}
		&=& \tr(B_1^TB_2B_2^TB_1) + \tr(B_1B_2^TB_2 B_1^T) -2\tr(B_1^TB_2B_1^TB_2).
\end{eqnarray}

\paragraph{Stiefel sectional curvature}
The Stiefel manifold of orthonormal $p$-frames can be realized as 
the quotient space $St(n, p) = SO(n)/SO(n-p)$.
The canonical projection and the left cosets are 
\begin{eqnarray*}
	\label{eq:Stiefel_quotient1}
	\Pi^{St}&:&SO(n)\to St(n,p)\\ 
	\label{eq:Stiefel_quotient_cosets} 
	[Q]=: \Pi^{St}(Q)  &=& \{Q\begin{pmatrix}
		I& 0\\
		0 & R
	\end{pmatrix}\mid R\in SO(n-p)\}.
\end{eqnarray*}
By splitting $Q = \begin{pmatrix}[c|c]
	U & U_\bot
\end{pmatrix}$,   $[Q]\in St(n,p)$ is uniquely determined by $U\in \R^{n\times p}$.\\
When lifting $[Q]\in St(n,p)$ to $Q\in SO(n)$, the vertical space is represented by
\begin{equation*}
	V_Q^{St} = 
	\bigl\{ \bar X\in \R^{n\times n}\mid \bar X = Q
	\begin{pmatrix}
		0 & 0\\
		0 & C
	\end{pmatrix},\\
	C\in \Skew(n-p)
	\bigr\} = Q\mathfrak{h}^{St}.
\end{equation*}
The associated horizontal space is
\begin{equation*}
	T_{[Q]}St(n,p) \cong 	H_Q^{St} =
	\bigl\{\bar X = Q
	\begin{pmatrix}
		A & -B^T\\
		B &  0
	\end{pmatrix} \mid \\
	A\in \Skew(p), B\in \R^{(n-p)\times p}
	\bigr\}= Q\mathfrak{m}^{St}.
\end{equation*}
For tangent vectors $X =\begin{pmatrix}
	A_1 & -B^T_1\\
	B_1 &  0
\end{pmatrix},
Y=			\begin{pmatrix}
	A_2 & -B^T_2\\
	B_2 &  0
\end{pmatrix}\in \mathfrak m^{St} \cong T_{[I]}St(n,p)$
that are orthonormal w.r.t. \eqref{eq:innerSOn}, 
a straightforward evaluation of \eqref{eq:curv_basis} yields
\begin{eqnarray}
	\nonumber
	K^{St}_c(X,Y) &=& 	\frac14 \| [ X, Y]_{{\mathfrak m}^{St}}  \|^2_{[I]} + \| [ X, Y]_{\mathfrak{h}^{St}} \|^2_{[I]}\\
		\nonumber
	&=&  \frac12\|B_2B_1^T - B_1B_2^T\|_F^2 +  \frac14\|B_1A_2 - B_2A_1\|_F^2\\  
	\label{eq:seccurv_Stief}
	& & + \frac18 \| [A_1,A_2] - (B_1^TB_2 - B_2^TB_1)\|_F^2.
\end{eqnarray}
This is in line with the result of \cite[Prop 4.2, eq. (34)]{nguyen2022curvature}.
%
%
\section{Matrix norm inequalities and curvature estimates}
\label{sec:curv_estimates}
In this section, we investigate the extremal behavior of the sectional curvature on the Grassmann and the Stiefel manifold.
It is known since \cite{Wong1968}, that the sectional curvature on the Grassmann manifold ranges in the interval $[0,2]$ with both the lower and the 
upper bound attained.
The upper bound can be established by using the matrix inequality 
\begin{equation}
	\label{eq:WuChen}
	\|B_1B_2^T - B_2B_1^T\|_F^2 \leq 2\|B_1\|_F^2\|B_2\|^2_F
\end{equation}
of Wu and Chen, \cite{WuChen1988}. 

Applying \eqref{eq:WuChen} to the terms \eqref{eq:seccurv_Grassmann1} and keeping $\|B_1\|_F=1=\|B_2\|_F$ in mind immediately gives
$0\leq K^{Gr}(X,Y)\leq 2$.
A key idea in \cite{WuChen1988} is to exploit the skew-symmetry of the real Schur form of $B_1B_2^T - B_2B_1^T$.
They also show by algebraic means that the inequality is sharp if and only if 
\[
	B_1 = \begin{pmatrix}[c|c]
		\begin{matrix} 0 & 1\\ \lambda & \mu \end{matrix} & \mathbf{0} \\
		\hline
		\mathbf{0}& \mathbf{0} \\ 
	\end{pmatrix},
	\quad
		B_2 = \begin{pmatrix}[c|c]
		\begin{matrix} \lambda & \mu\\ 0 & -1 \end{matrix} & \mathbf{0} \\
		\hline
		\mathbf{0}&\mathbf{0} \\ 
	\end{pmatrix},
\]
up to an orthogonal transformation and scaling.
Hence, under the normalization constraint $\|B_1\|_F=1=\|B_2\|_F$, both terms in \eqref{eq:seccurv_Grassmann1} attain their upper bound of $2$ simultaneously for the rank-2 matrices 
\[
		B_1 = \frac{1}{\sqrt{2}}\begin{pmatrix}[c|c]
		\begin{matrix} 0 & 1\\ 1 & 0 \end{matrix} & \mathbf{0} \\
		\hline
		\mathbf{0}& \mathbf{0} \\ 
	\end{pmatrix},
	\quad
	B_2 = \frac{1}{\sqrt{2}}\begin{pmatrix}[c|c]
		\begin{matrix} 1 & 0\\ 0 & -1 \end{matrix} & \mathbf{0} \\
		\hline
		\mathbf{0}&\mathbf{0} \\ 
	\end{pmatrix}.
\]
Observe that $\langle B_1,B_2\rangle_F=0$. Hence, this matrix pair forms an ONB and so do the associated Grassmann tangent vectors.

To gain more insight on how the various trace terms contribute to the overall curvature value, 
we will re-establish this result via an optimization approach.
To this end, we will consider the three trace terms 
\[
\tr(B_1^TB_2B_2^TB_1), \mbox{ }	\tr(B_1B_2^TB_2B_1^T),  \mbox{ and } -2\tr(B_1^TB_2B_1^TB_2) \quad B_1, B_2\in \R^{(n-p)\times p},
\]
in the curvature expression \eqref{eq:seccurv_Grassmann2} 
separately. Eventually, this will leads to an alternative proof and a refinement of the matrix inequality of Wu and Chen, \cite{WuChen1988}.
We start with a preparatory lemma that improves on the classical submultiplicativity property $\|AB\|_F\leq \|A\|_F\|B\|_F$.
\begin{lemma}
	\label{lem:2norm_Fnorm}
 Let $A\in \R^{n\times m}$,  $B\in \R^{m\times p}$. 
	Then
	\begin{equation}
    \label{eq:2norm_Fnorm}
	    		\|AB\|_F \leq \min \{\|A\|_2 \|B\|_F, \|A\|_F \|B\|_2\}.
	\end{equation}

 If either $A$ or $B$ is skew-symmetric, then
 \begin{equation}
     \label{eq:2norm_Fnorm_skew}
     \|AB\|_F \leq \frac{1}{\sqrt{2}} \|A\|_F \|B\|_F.
 \end{equation}

\end{lemma}
\begin{proof}
	On \eqref{eq:2norm_Fnorm}: First, consider $D\in\R^{m\times m}$ diagonal.
	Write $B=(b_1,\ldots,b_p)$ column-wise. It holds 
	\begin{eqnarray*}
		\|DB\|^2_F&=&\tr(B^TD^2B)=\sum_{j=1}^p b_j^TD^2b_j=\sum_{j=1}^p \|b_j\|_2^2 \frac{b_j^TD^2b_j}{b_j^Tb_j}\\
		&\leq& \max_j\{d_j^2\} \sum_{j=1}^p \|b_j\|_2^2 =\|D\|_2^2 \|B\|_F^2.
	\end{eqnarray*}
	The general case can be reduced to this case.
	 Let $U\Sigma V^T=A$ be the full SVD of $A$.
	\begin{eqnarray*}
	\|AB\|_F^2	&=& \|U\Sigma V^T B\|_F^2 = \tr(B^TV\Sigma^T\Sigma V^TB)\\
	&=& \|\sqrt{\Sigma^T\Sigma} V^TB\|_F^2 \stackrel{D= \sqrt{\Sigma^T\Sigma}}{=} \sigma_1^2 \|V^TB\|_F^2 = \|A\|_2^2 \|B\|_F^2.
	\end{eqnarray*}
	When working with the SVD of $B$, the same argument may be applied to $\|B^TA^T\|^2_F$ 
	and yields $\|B^TA^T\|^2_F\leq  \|B^T\|_2^2 \|A^T\|_F^2 =  \|B\|_2^2 \|A\|_F^2$.\\
    The second inequality \eqref{eq:2norm_Fnorm_skew} comes as a corollary. W.l.o.g. assume that $A\in\Skew(m)$. Then, all eigenvalues of $A$ are either zero or form purely imaginary, complex conjugate pairs. The singular values of $A$ are the absolute values of the eigenvalues of $A$. Therefore, the singular values also come in pairs. Let $\sigma_1 = \sigma_2\geq\dots\geq\sigma_{r-1} = \sigma_r>0$ be the non-zero singular values of $A$. We obtain
    \[\|A\|_2^2 = \sigma_1^2 = \frac12\left(\sigma_1^2 + \sigma_2^2\right)\leq\frac12\left(\sigma_1^2 + \sigma_2^2 + \dots + \sigma_r^2\right) = \frac12\|A\|_F^2.\]
    Hence, \eqref{eq:2norm_Fnorm} yields $\|AB\|_F \leq \|A\|_2 \|B\|_F \leq \frac{1}{\sqrt{2}} \|A\|_F \|B\|_F$.    
\end{proof}
Alternatively, \eqref{eq:2norm_Fnorm} can also be established as a consequence of \cite[Lemma 1]{WangKuoHsu1986}.

The following lemma is obvious.
\begin{lemma}
	\label{lem:trace_max}
	Let $m\geq p$ and consider $B_2\in \R^{m\times p}$ as fixed. Let $U\Sigma V^T=B_2$ be the (full) SVD of $B_2$ with 
	$\Sigma = \begin{pmatrix}
		\Sigma_p\\
		\mathbf{0}
	\end{pmatrix}\in\R^{m\times p}$ and $\Sigma_p=\mbox{diag}(\sigma_1,\ldots,\sigma_p)$.
	\begin{enumerate}[(1.)]
		\item The global optimum of
		\[
			\max_{B_1\in\R^{m\times p}} \tr(B^T_1B_2B_2^TB_1) \quad \mbox{s.t. } \|B_1\|_F=1
		\]
		is $\sigma_1^2$ and is attained for any normalized $B_1$ such that $\tilde B_1 := U^T B_1 V$
		only features scaled copies of the first unit vector $e_1 = (\pm1,0,\ldots,0)^T \in \R^{m}$ as columns.
		\item The global optimum of
		\[
		\max_{B_1\in\R^{m\times p}} \tr(B_1B_2^TB_2B^T_1) \quad \mbox{s.t. } \|B_1\|_F=1
		\]
		is $\sigma_1^2$ and is attained for any normalized $B_1$ such that $\tilde B_1 := U^T B_1 V$
		only features scaled copies of the first unit vector $e_1 = (\pm1,0,\ldots,0)^T\in \R^{p}$ as rows.
	\end{enumerate}
\end{lemma}
The next lemma concerns the third trace term in \eqref{eq:seccurv_Grassmann2}.
\begin{lemma}
	\label{lem:B1B2B1B2_term}
	Let $m\geq p$ and let $B_2 = U\Sigma V^T\in \R^{m\times p}$ (not necessarily normalized). 
	Then
	\begin{eqnarray*}
		\max_{B_1\in\R^{m\times p}, \|B_1\|_F=1} \tr(B^T_1B_2 B^T_1B_2) &=& \sigma_{1}^2,
		\\
		\min_{B_1\in\R^{m\times p}, \|B_1\|_F=1} \tr(B^T_1B_2 B^T_1B_2) &=& -\sigma_1\sigma_2.
	\end{eqnarray*}
	The global maximum and minimum are attained for $B_1^{+} = U\tilde B_1^{+} V^T$ and $B_1^{-} = U\tilde B_1^{-} V^T$, where
    \begin{equation}
    	\label{eq:opt_B1B2B1B2_term}
    \tilde B_1^{+} = \pm \begin{pmatrix}[c|c]
		\begin{matrix}  1 & 0\\ 0 & 0 \end{matrix} & \mathbf{0} \\
		\hline
		\mathbf{0}& \mathbf{0} \\ 
	\end{pmatrix}, \mbox{ and }
	\tilde B_1^{-} =\pm \frac{ 1}{\sqrt{2}}\begin{pmatrix}[c|c]
		\begin{matrix} 0 & 1\\ -1 & 0 \end{matrix} & \mathbf{0} \\
		\hline
		\mathbf{0}& \mathbf{0} \\ 
	\end{pmatrix}, \mbox{ respectively.}
	\end{equation}
	For $B_1, B_2$ both of unit Frobenius norm, 
	\[
		-\frac12 \leq \tr(B_1^TB_2 B^T_1B_2) \leq 1.
	\]
	In this case, the global extrema are attained for $B_1^+,B_1^-$ as above and $B_2^+$ of rank one ($\sigma_1=1$) 
	and $B_2^-$ of rank 2 ($\sigma_1 = \sigma_2 = \frac{1}{\sqrt{2}}$).
\end{lemma}
\begin{remark}
	The extrema may not be isolated nor are they necessarily unique (up to the trace-preserving transformation).
	If $m=p$, $B_2 = \sigma I$, then any normalized symmetric $B_1$ yields the maximum $\tr(B^T_1B_2 B^T_1B_2)= \sigma^2\|B_1\|^2_F  = \sigma^2_1$.
	Likewise, any normalized skew-symmetric  $B_1$ yields the minimum $\tr(B^T_1B_2 B^T_1B_2)= -\sigma^2\|B_1\|^2_F  = -\sigma^2=(-\sigma_1\sigma_2)$.
	However, for both $B_1$ and $B_2$ of unit Frobenius norm, the global extrema are attained only for the matrices in \eqref{eq:opt_B1B2B1B2_term}
	(up to the trace preserving transformation).
\end{remark}
\begin{proof}[\Cref{lem:B1B2B1B2_term}]
	Let $B_2 = U\Sigma V^T\in \R^{m\times p}$ be the full SVD of $B_2$. Under the norm-preserving bijection $B\mapsto U^TBV=:\tilde B$, 
	the trace optimization problem $\min_{B\in\R^{m\times p}, \|B\|_F=1} \pm \tr(B^TB_2B^TB_2)$ becomes
	 $\min_{\tilde B\in\R^{m\times p}, \|\tilde B\|_F=1} \pm\tr(\tilde B^T\Sigma\tilde B^T\Sigma)$.

	Write $\Sigma= \begin{pmatrix}
		\Sigma_p\\
		\mathbf{0}
	\end{pmatrix}$ and 
    $\tilde B=\begin{pmatrix}
    	\tilde B_p\\
    	\tilde B_{m-p}
    	\end{pmatrix}$, with empty lower blocks if $m=p$.
    As a preliminary, observe that $\tr(\tilde B_p^T \tilde B_p^T)= \langle \tilde B_p, \tilde B^T_p\rangle_F = \sum_{j=1}^p \sum_{k=1}^p b_{jk}b_{kj} $.
    It holds
	\begin{align}
		\label{eq:trace_B1B2B1B2_1}
		\tr(\tilde B^T\Sigma\tilde B^T\Sigma) &= \sum_{j=1}^p \sum_{k=1}^p \sigma_j\sigma_k b_{jk}b_{kj} \\
		\nonumber 
		&= \sigma_1^2b_{11}^2 + \sum_{k=2}^p \sigma_1\sigma_k b_{1k}b_{k1} +\sum_{j=2}^p\left( \sigma_j\sigma_1 b_{j1}b_{1j} + \sum_{k=2}^p \sigma_j\sigma_k b_{jk}b_{kj}\right)\\
		\nonumber
		&\leq \sigma_{1}^2 \langle \tilde B_p, \tilde B_p^T\rangle_F
		\leq \sigma_{1}^2 \|\tilde B\|_F^2 = \sigma_{1}^2.
	\end{align}
	The maximum is attained, if all weight in $\tilde B$ is put on the upper diagonal entry, i.e., for $\tilde B^{+}$ as in the statement of the lemma.
	If $\sigma_1>\sigma_2$, then the maximum is isolated.
	
	Now on the minimum. 
	Reconsider \eqref{eq:trace_B1B2B1B2_1} and extract the `diagonal terms', 
		\[
		\tr(\tilde B^T\Sigma\tilde B^T\Sigma) = \sum_{j=1}^p \sum_{k=1,k\neq j}^p \sigma_j\sigma_k b_{jk}b_{kj} + \sum_{j=1}^p{\sigma_j^2 b_{jj}^2}.
		\]
	One sees that the diagonal terms in $\tilde B_p$ only make nonnegative contributions to the trace total.
	All remaining terms become non-positive, if $b_{jk}$ and $b_{kj}$ are of opposite sign for all $k\neq j$.
	In this case and with $b_{jj}=0$, it holds
		\begin{eqnarray}
		\nonumber
		\tr(\tilde B^T\Sigma\tilde B^T\Sigma) =& -\sum_{j=1}^p \sum_{k=1,k\neq j}^p \sigma_j\sigma_k b_{jk}b_{kj} 
		&\geq -\sigma_1\sigma_2 \sum_{j=1}^p \sum_{k=1,k\neq j}^p b_{jk}b_{kj}\\
		\nonumber
		\geq& -\sigma_1\sigma_2 \langle \tilde B_p, \tilde B_p^T\rangle_F
		&\geq -\sigma_1\sigma_2 \|\tilde B\|_F^2 = -\sigma_{1}\sigma_2.
	\end{eqnarray}
	The pairing $\sigma_{1}\sigma_2$ of the largest singular values features only  as a factor in front of the product $b_{12}b_{21}$.
	Therefore, the estimate is sharp if all weight in $\tilde B$ is placed on these terms.
	Consider $\tilde B = \begin{pmatrix}[c|c]
		\begin{matrix} 0 & b_{12}\\ b_{21} & 0 \end{matrix} & \mathbf{0} \\
		\hline
		\mathbf{0}& \mathbf{0} \\ 
	\end{pmatrix}$ with $b_{12}^2 +b_{21}^2 =1$.  
	This yields
	\begin{eqnarray*}
			\tr(\tilde B^T\Sigma\tilde B^T\Sigma) &= & \tr(\begin{pmatrix} 0 & b_{12}\\ b_{21} & 0 \end{pmatrix}
			                                               \begin{pmatrix} \sigma_1 & 0\\ 0 & \sigma_2 \end{pmatrix}
			                                               \begin{pmatrix} 0 & b_{12}\\ b_{21} & 0 \end{pmatrix}
			                                               \begin{pmatrix} \sigma_1 & 0\\ 0 & \sigma_2 \end{pmatrix})\\
			                                       &= &2\sigma_1\sigma_2b_{12}b_{21}.
	\end{eqnarray*}
	The product $b_{12}b_{21}$ gets extremal for $b_{12},b_{21} = \pm\frac{1}{\sqrt{2}}$, the associated trace minimum is
	$-\sigma_1\sigma_2$ and is attained when $b_{12}$ and $b_{21}$ feature opposite signs.
	The global minimum is isolated, if $\sigma_2>\sigma_3$.
	
	Under the additional normalization of $\|B_2\|_F=1$, the global minimum of $-\frac12$ is attained, if $\sigma_1=\sigma_2 = \frac{1}{\sqrt{2}}$
	(which enforces $\sigma_3,\ldots,\sigma_p =0$).
	The global maximum of a value of $1$ is attained if $\sigma_1= 1$ (which enforces $\sigma_2,\ldots,\sigma_p =0$).
\end{proof}

The next theorem includes (and refines) the inequality of Wu and Chen. The proof is different and is based on an optimization approach.
\begin{theorem}
	\label{thm:WuChen}
	For $B_1,B_2\in \R^{m\times p}$, with SVDs 
	$B_1 = U_1 P V_1^T\in \R^{m\times p}$, 
	$B_2 = U_2\Sigma V_2^T\in \R^{m\times p}$, 
	where the upper left diagonal blocks of the singualar values matrices $P,\Sigma$ are
	$P_p= \operatorname{diag}(\rho_1,\ldots,\rho_p)$, and $\Sigma_p= \operatorname{diag}(\sigma_1,\ldots,\sigma_p)$,
	respectively,
	it holds
	\begin{align}
	\label{eq:WuChenimproved}
			 \tr(B_1^TB_2B_2^TB_1) - \tr(B_1^TB_2B_1^TB_2) 
		&\leq \min\{ \|B_1\|_F^2 (\sigma_1^2 + \sigma_2^2), \|B_2\|_F^2 (\rho_1^2 + \rho_2^2) \},\\
	\label{eq:WuChenimproved1}
	\tr(B_1B_2^TB_2B_1^T) - \tr(B_1^TB_2B_1^TB_2) 
	&\leq \min\{ \|B_1\|_F^2 (\sigma_1^2 + \sigma_2^2), \|B_2\|_F^2 (\rho_1^2 + \rho_2^2) \}	.
	\end{align}
	As a consequence,
	\begin{equation}
		\label{eq:WuChen_revisited}
		\frac12 \|B_1^TB_2 - B_2^TB_1\|_F^2 \leq \|B_1\|_F^2\|B_2\|_F^2, \quad 
		\frac12 \|B_1B_2^T - B_2B_1^T\|_F^2 \leq \|B_1\|_F^2\|B_2\|_F^2.
	\end{equation}
\end{theorem}
Note that \eqref{eq:WuChen_revisited} is the Wu-Chen matrix inequality of \cite{WuChen1988}.
While \eqref{eq:WuChenimproved}, \eqref{eq:WuChenimproved1} do not give tighter general bounds, they can be significantly tighter in special situations.
For example, if $B_2\in St(m,p)$ is column-orthogonal, then all singular values of $B_2$ are $1$ and  $\|B_2\|_F^2 = \tr(I)= m$. 
The bound of \eqref{eq:WuChenimproved} is 
$\|B_1^TB_2 - B_2^TB_1\|_F \leq \sqrt{2}\sqrt{2}\|B_1\|_F$, while that of
\eqref{eq:WuChen_revisited} gives $\|B_1^TB_2 - B_2^TB_1\|_F \leq  \sqrt{2}\sqrt{m}\|B_1\|_F$, which grows with the dimension $m$.
\begin{proof}
	It holds $\frac12 \|B_1^TB_2 - B_2^TB_1\|_F^2 = \tr(B_1^TB_2B_2^TB_1) - \tr(B_1^TB_2B_1^TB_2) $
	and $\frac12 \|B_1B_2^T - B_2B_1^T\|_F^2 = \tr(B_1B_2^TB_2B_1^T) - \tr(B_1^TB_2B_1^TB_2) $.
	We normalize and 
	apply again the coordinate change based on the SVD data of $B_2$ to obtain
	\begin{align*}
		&\tr(B_1^TB_2B_2^TB_1) - \tr(B_1^TB_2B_1^TB_2)\\
		=&\|B_1\|_F^2\tr\left(\frac{B_1^T}{\|B_1\|_F}B_2B_2^T\frac{B_1}{\|B_1\|_F}\right)
		 - \|B_1\|_F^2 \tr\left(\frac{B_1^T}{\|B_1\|_F}B_2 \frac{B_1^T}{\|B_1\|_F}B_2\right)\\
		=&\|B_1\|_F^2 \left(\tr(\tilde B^T\Sigma\Sigma^T\tilde B) - \tr(\tilde B^T\Sigma \tilde B^T\Sigma)\right), \quad 
		\tilde B = U^T(B_1/\|B_1\|_F)V, \hspace{0.2cm}\|\tilde B\|_F=1.
	\end{align*}
	By the preparatory lemmata, it is clear that all weight in the (normalised) matrix $\tilde B$ must be concentrated on the upper $(2\times 2)$-diagonal block.
 	Therefore, we look for $\tilde B = \begin{pmatrix}
 		a & b\\
 		c & d
 	\end{pmatrix}$
 	with the optimal balance between the parameters $a,b,c,d$ that maximizes the combination of trace terms.
 	Yet, a quick calculation shows that the diagonal entries $a,d$ do not contribute to the result.
 	We have
 	\[
 	 \tr(\tilde B^T\Sigma\Sigma^T\tilde B) - \tr(\tilde B^T\Sigma \tilde B^T\Sigma)
 	 = \sigma_1^2b^2 + \sigma_2^2 c^2 -2\sigma_1\sigma_2bc = (\sigma_1b - \sigma_2 c)^2.
 	\]
 	To maximize this term, all weight in $\tilde B$ must be on the off-diagonal terms so that 
 	$a,d=0$ and $b^2 +c^2=1$. Let $S^1 = \{(b,c)\in\R^2\mid b^2 + c^2 =1\}$ be the unit circle, parameterized by
 	$\gamma:t\mapsto (b(t),c(t)= (\cos(t), \sin(t))$.
 	The function $f:[-\pi,\pi) \to \R, t\mapsto (\sigma_1 \cos(t) - \sigma_2 \sin(t))^2$ attains
	its global maximum on $S^1$ at $t_* = \arctan\left(-\frac{\sigma_2}{\sigma_1}\right)$.
	This yields $b_* = \frac{\sigma_1}{\sqrt{\sigma_1^2 +\sigma_2^2}}, c_* = \frac{-\sigma_2}{\sqrt{\sigma_1^2 +\sigma_2^2}}$
	and a global maximum of 
	\begin{equation}
        \label{eq:max_WuChen14}
		(\sigma_1b_* - \sigma_2 c_*)^2 = (\sigma_1^2 + \sigma_2^2).
	\end{equation}
	Combined, it holds
	\[
		\tr(B_1^TB_2B_2^TB_1) - \tr(B_1^TB_2B_1^TB_2) \leq  \|B_1\|_F^2(\sigma_1^2 + \sigma_2^2)\leq \|B_1\|_F^2\|B_2\|_F^2.
	\]
	The roles of $B_1$ and $B_2$ can be exchanged.
	The same reasoning applies to \eqref{eq:WuChenimproved1} but with the transpose of $B_1,B_2$.
	This yields
	\begin{equation}
        \label{eq:max_WuChen15}
		 \tr(\tilde B\Sigma^T\Sigma\tilde B^T) - \tr(\tilde B^T\Sigma \tilde B^T\Sigma)
		=  (\sigma_2b - \sigma_1 c)^2.
	\end{equation}

	The global optimum is the same value of $(\sigma_1^2 + \sigma_2^2)$, but is attained for
	 $b_* = \frac{\sigma_2}{\sqrt{\sigma_1^2 +\sigma_2^2}}$, $c_* = \frac{-\sigma_1}{\sqrt{\sigma_1^2 +\sigma_2^2}}$,
	 which corresponds to the transpose of the maximizer of \eqref{eq:WuChenimproved}.
	 Both inequalities become simultaneously sharp for the same input pair $\tilde B, \Sigma$, if $\sigma_1=\sigma_2$.
\end{proof}

\subsection{The classical Grassmann sectional curvature bounds}
With \Cref{thm:WuChen}, the global bound of $K^{Gr}\leq 2$ of the sectional curvature on $Gr(n,p)$ is an immediate consequence of
\eqref{eq:seccurv_Grassmann1}, \eqref{eq:seccurv_Grassmann2} (as was already clear from \cite{WuChen1988}, and clear to Wong in 1968, \cite{Wong1968}).
The curvature formulas \eqref{eq:seccurv_Grassmann1}, \eqref{eq:seccurv_Grassmann2} hold for $\{B_1,B_2\}$ forming an ONB.
With $B_2 = U\Sigma V^T$, under the transformation $B_1\mapsto U^TB_1V=:\tilde B_1$, this is equivalent to
$\{\tilde B_1, \Sigma\}$ being an ONB.
According to \Cref{thm:WuChen}, the matrices, for which both matrix inequalities are simultaneously sharp are
\begin{equation}
\label{eq:curv_maximizers}
	\tilde B_1 = \pm \begin{pmatrix}[c|c]
		\begin{matrix}  0 & \frac{1}{\sqrt{2}}\\ \frac{-1}{\sqrt{2}} & 0 \end{matrix} & \mathbf{0} \\
		\hline
		\mathbf{0}& \mathbf{0} \\ 
	\end{pmatrix}, \mbox{ and }
	\tilde \Sigma = \begin{pmatrix}[c|c]
		\begin{matrix} \sigma & 0\\ 0 & \sigma \end{matrix} & \mathbf{0} \\
		\hline
		\mathbf{0}& \mathbf{0} \\ 
	\end{pmatrix}.
\end{equation}
This pair of matrices is orthogonal; it is an ONB if $\sigma = \frac{1}{\sqrt{2}}$.
Observe that $\tilde B_1^T \Sigma$ is skew-symmetric.
The tangent plane spanned by the associated tangent vectors is the only tangent plane of maximum sectional curvature.
\subsection{Bounds on the canonical sectional curvature on Stiefel}
In this section, we establish a global upper bound on the sectional curvature on $St(n,p)$ under the canonical metric for all Stiefel manifolds $n\geq p$ that are at least two-dimensional. (Otherwise, the concept of sectional curvature does not apply.)
We start from \eqref{eq:seccurv_Stief} for an orthonormal pair of Stiefel tangents
\[
	X=\begin{pmatrix}
		A_1 & -B_1^T\\
		B_1 & \mathbf{0}
	\end{pmatrix},\quad 
	Y=
	\begin{pmatrix}
		A_2 & -B_2^T\\
		B_2 & \mathbf{0}
	\end{pmatrix} \in T_{[I]}St(n,p).
\]
Orthonormality in the canononical Stiefel metric means that
\[
	\frac12 \|A_j\|_F^2 + \|B_j\|^2_F = 1, \hspace{0.1cm}j=1,2,\quad \frac12 \tr(A_1^TA_2) + \tr(B_1^TB_2) = 0.
\]
From \cite[Lemma 2.5]{Ge2014}, we have norm bounds for the commutator bracket of skew-symmetric matrices, 
\begin{eqnarray}
\label{eq:comm_skew4}
	\|[A_1,A_2]\|_F^2\leq& \|A_1\|_F^2\|A_2\|_F^2, \quad&\forall p\geq 4,\\
 \label{eq:comm_skew3}
	\|[A_1,A_2]\|_F^2\leq& \frac12\|A_1\|_F^2\|A_2\|_F^2,\quad& p=3,\\
 \label{eq:comm_skew2}
	\|[A_1,A_2]\|_F^2=& 0, \quad& p=2.
\end{eqnarray}
The last one is obvious, because $(2\times 2)$-skew symmetric matrices necessarily commute; \eqref{eq:comm_skew3} is a straightforward consequence of the interplay between skew-symmetric $(3\times 3)$-matrices and their representation with $3$-vectors, see \Cref{app:curv_special_dim} for a short recap.

Write $\alpha_1 = \|A_1\|_F, \alpha_2= \|A_2\|_F,\beta_1 = \|B_1\|_F, \beta_2= \|B_2\|_F$ and note that for normalized tangent vectors
$\beta_j^2 = 1-\frac12\alpha_j^2$, $\beta_j\in[0,1], \alpha_j\in [0,\sqrt{2}]$.

Because $A_1,A_2$ are skew-symmetric, eq. \eqref{eq:2norm_Fnorm_skew} of \Cref{lem:2norm_Fnorm} yields  the following refined estimate for the `one-fourth'-term in \eqref{eq:seccurv_Stief}
\begin{eqnarray*}
\frac14\|B_1A_2 - B_2A_1\|_F^2 &\leq& \frac14\left(\frac{1}{\sqrt{2}}\|B_1\|_F\|A_2\|_F + \frac{1}{\sqrt{2}}\|B_2\|_F\|A_1\|_F\right)^2\\
&=& \frac18\left(\beta_1^2\alpha_2^2 + \beta_2^2\alpha_1^2 + 2\alpha_1\alpha_2\beta_2\beta_2\right).
\end{eqnarray*}
With the above inequalities and \Cref{thm:WuChen}, a global bound for the Stiefel curvature is
\begin{eqnarray}
	\nonumber
	K^{St}_c(X,Y) 
	&=&  \frac12\|B_2B_1^T - B_1B_2^T\|_F^2 +  \frac14\|B_1A_2 - B_2A_1\|_F^2\\  
 	\nonumber
	& & + \frac18 \| [A_1,A_2] - (B_1^TB_2 - B_2^TB_1)\|_F^2\\
    \label{eq:brute_force_ineq}
	&\leq &
	\frac54 + \frac{5}{16}\alpha_1^2\alpha_2^2 - \frac12 (\alpha_1^2 + \alpha_2^2)
	+\frac{1+\sqrt{2}}{4}\alpha_1\alpha_2\sqrt{1-\frac12\alpha_1^2}\sqrt{1-\frac12\alpha_2^2}.
\end{eqnarray}
This function in $(\alpha_1,\alpha_2)$ has an isolated local maximum at $\alpha_1=0=\alpha_2$ with a corresponding function value of $\frac54$, which is the global maximum in the 
admissible range of $(\alpha_1,\alpha_2)\in[0,\sqrt{2}]^2$. For a verification see \Cref{app:auxresults}.
For most Stiefel manifolds, more precisely, for all Stiefel manifolds with $p\geq 2$, $n\geq p+2$, this bound is tight. 
This main result is detailed in the next theorem.
For the sake of completeness, the remaining cases
are also included, even though they are not new. The cases $p=2,n=3$ and $p=2, n>3$ are explicitly treated in \cite[Prop. 6.1]{nguyen2022curvature} for a parametric family of metrics \cite{HueperMarkinaLeite2020, ZimmermannHueper2022} that include the canonical metric as a special case.
Since $St(n,1)\cong S^{n-1}$, $St(n,n-1)\cong O(n)$ and $St(n,n)\cong SO(n)$, one may argue that $St(n,p)$, $p\geq 2$, $n\geq p+2$ are the only `true' Stiefel manifolds.
\begin{theorem}
	\label{thm:StiefelCurvBound}
	The sectional curvature under the canonical metric on the Stiefel manifold $St(n,p)$, $n\geq p$, is globally bounded by 
	\[
	0\leq K^{St}_c(X,Y) \leq \frac54.
	\]
 \begin{enumerate}
     \item For $p\geq 2, n\geq p+2$, the bound is sharp.
 Up to trace-preserving transformations, the maximum curvature is attained only for the tangent plane spanned by $\left\{X= \begin{pmatrix}
	0 & -B_1^T\\
	B_1 & 0
\end{pmatrix},Y=\begin{pmatrix}
	0 & -B_2^T\\
	B_2 & 0
\end{pmatrix}\right\}$,
where 
\begin{equation}
\label{eq:secmaximizers_n4p2}
	B_1 = \pm\frac{1}{\sqrt{2}} \begin{pmatrix}[c|c]
		\begin{matrix}  0 & 1\\ -1 & 0 \end{matrix} & \mathbf{0} \\
		\hline
		\mathbf{0}& \mathbf{0} \\ 
	\end{pmatrix}, \mbox{ and }
	B_2 = \pm\frac{1}{\sqrt{2}}\begin{pmatrix}[c|c]
		\begin{matrix} 1 & 0\\ 0 & 1 \end{matrix} & \mathbf{0} \\
		\hline
		\mathbf{0}& \mathbf{0} \\ 
	\end{pmatrix}.
\end{equation}
    \item 	If $\rank(B_1) = 1$ or $\rank(B_2)= 1$ holds for the $B$-blocks in $X$ or $Y$, respectively, then 
	\[
	K^{St}_c(X,Y) \leq 1.
	\]
    The bound is attained for the matrices stated in eq. \eqref{eq:rank1_bound_argmax} below.
    \item For $p=1$ and $n\geq 3$, 
    $$K^{St}_c(X,Y) \leq 1$$ 
    and the bound is sharp. It is attained for the matrices specified in eq. \eqref{eq:rank1_bound_argmax} if these are reduced to their first column.
    \item For $(n,p)=(3,2)$ or $(n,p)=(3,3)$, the sectional curvature has a constant value of $$K^{St}_c(X,Y)\equiv\frac14.$$
    \item For $n\geq 4$ and  $p=n-1$ or $p=n$,
    $$K^{St}_c(X,Y)\leq\frac12$$
    and the bound is sharp.
 \end{enumerate}
 For $n=2,p=2$ or $n=2,p=1$, $St(n,p)$ is one-dimensional so that the concept of sectional curvatures does not apply.
\end{theorem}
\begin{proof}
    On 1.: The global bound of $\frac54$ is already established by the maximum of \eqref{eq:brute_force_ineq}. A direct calculation shows that the bound is attained for the tangent plane associated with the matrices in \eqref{eq:secmaximizers_n4p2}. Hence, the bound is sharp for all Stiefel manifolds that fit these matrices dimension-wise.
    Moreover, the matrices in \eqref{eq:secmaximizers_n4p2} are the only matrices (up to trace-preserving transformations) for which the bound is achieved. This can be established analogously to the Grassmann manifold case, since the maximum curvature is attained for zero skew-symmetric blocks $A_j=0$.\\
    On 2.: If one of the matrices $B_1$ or $B_2$ is of rank one, 
	then from the proof of \Cref{thm:WuChen} one can deduce that
	\begin{eqnarray*}
		\frac12\|B_2B_1^T - B_1B_2^T\|_F^2 + \frac12 \|B_1^TB_2 - B_2^TB_1\|_F^2
		\leq \|B_1\|_F^2\|B_2\|_F^2.
	\end{eqnarray*}
  Compared to the general case, the bound is improved by a factor of $2$. As before, the maximum curvature is  attained for zero skew-symmetric blocks $A_j=0$.
  Hence, the curvature is bounded by 
  \begin{eqnarray*}
  K^{St}_c(X,Y) &\leq&   \frac12\|B_2B_1^T - B_1B_2^T\|_F^2  
	                + \frac18 \|(B_1^TB_2 - B_2^TB_1)\|_F^2\\
                &=&  \frac38\|B_2B_1^T - B_1B_2^T\|_F^2  \\
                &&
	                + \frac14 \left(
                 \frac12\|B_2B_1^T - B_1B_2^T\|_F^2 + 
                 \frac12\|(B_1^TB_2 - B_2^TB_1)\|_F^2\right)\\
                &\leq& \frac34 \|B_1\|_F^2\|B_2\|_F^2 + \frac14\|B_1\|_F^2\|B_2\|_F^2 \leq 1.
  \end{eqnarray*}
The rank-one bound is attained for
      \begin{equation}
    	\label{eq:rank1_bound_argmax}
    B_1= \pm \begin{pmatrix}[c|c]
		\begin{matrix}  0 & 0\\ 1 & 0 \end{matrix} & \mathbf{0} \\
		\hline
		\mathbf{0}& \mathbf{0} \\ 
	\end{pmatrix}, \mbox{ and }
	B_2 =\pm \begin{pmatrix}[c|c]
		\begin{matrix} 1 & 0\\ 0 & 0 \end{matrix} & \mathbf{0} \\
		\hline
		\mathbf{0}& \mathbf{0} \\ 
	\end{pmatrix}.
	\end{equation}
The matrix $B_1$ may feature multiple copies of the canonical unit vector $e_1^T$ as rows.\\
On 3.: This is a direct consequence of item 2.\\
On 4.: See \cite[Prop. 6.1]{nguyen2022curvature}
or see \Cref{app:curv_special_dim} for a direct proof.\\
On 5.: This is clear because $St(n,n-1)\cong O(n)$, $St(n,n)\cong SO(n)$.
    In both these cases, $[X,Y]_{\mathfrak{m}}=[X,Y], [X,Y]_{\mathfrak{h}}=0$ so that the curvature formula \eqref{eq:curv_basis} reduces to \eqref{eq:curv_biinv_Liegroup}. For $O(n)$ and $SO(n)$, the formula from \Cref{rem:curvSOn_fro} applies. The global bound of $\frac12$ stems from \eqref{eq:comm_skew4}.
\end{proof}

\begin{remark}
The bound of \Cref{thm:StiefelCurvBound} was already correctly guessed in \cite[p. 94, 95]{Rentmeesters2013}, but the derivation there is not correct. 
It was overlooked that tangent matrices $X$ of unit norm 
w.r.t. the Riemannian metric $1=\|X\|^2 = \frac12\|X\|_F^2$ have a
Frobenius norm of $\|X\|_F = \sqrt{2}$.
In \cite[Prop. 6.2]{nguyen2022curvature} it is shown that the sectional curvature range of the Stiefel manifold under the canonical metric {\em includes} the interval $[0,\frac54]$. But to argue that this is the exact range, the author resorted to the result of \cite[p. 94, 95]{Rentmeesters2013}, which was lacking a valid proof.
\end{remark}
\subsection{Bounds on the Euclidean sectional curvature on Stiefel}
Besides the canonical metric that follows from equipping the orthogonal group with the Euclidean metric $\langle X,Y\rangle_Q = \frac12\langle X,Y\rangle_F $ and the quotient space representation $St(n,p) = O(n)/O(n-p)$,
the most important metric for applications is obtained from considering the Stiefel manifold as an embedded submanifold of the Euclidean matrix space, $St(n,p)\subset \R^{n\times p}$. Tangent vectors are then rectangular matrices $\Delta = Q\begin{pmatrix}
    A\\
    B
\end{pmatrix} = UA+U_\bot B\in\R^{n\times p}$, where $Q=\begin{pmatrix}
    U& U_\bot
\end{pmatrix}\in O(n)$, $A\in \Skew(p), B\in \R^{(n-p)\times p}$.
The natural Riemannian metric in this setting is obtained from restricting the Euclidean inner product\footnote{Caution: Here {\em without} the factor $\frac12$!} to the Stiefel tangent spaces. This yields the so-called {\em Euclidean metric} on $St(n,p)\subset\R^{n\times p}$,
\[
\langle \Delta,\tilde \Delta\rangle^{St}_e
= \tr\left(\begin{pmatrix}
    A^T & B^T
\end{pmatrix}Q^T Q
\begin{pmatrix}
    \tilde A\\
    \tilde B
\end{pmatrix}\right)
= \tr(A^T\tilde A) + \tr(B^T\tilde B).
\]
The differences between this metric and the canonical one have been discussed in \cite[p. 313]{EdelmanAriasSmith1999}.
Most notably, this metric is not directly related to a quotient space construction so that the formulas of \Cref{sec:LieGroups} do not readily apply.
One option for calculating curvatures is via direct computations and the geometry of submanifolds, see \cite[Section 8]{Lee2018riemannian}.
Yet, the canonical and the Euclidean metric are included in a parametric family of metrics considered in \cite[Def. 3.1, Remark 1]{HueperMarkinaLeite2020}, where the authors show that all these metrics stem from a bi-invariant metric associated with a different quotient space representation, namely $St(n,p)\cong (O(n)\times O(p))/(O(n-p)\times O(p))$.
The associated sectional curvature is computed in \cite{nguyen2022curvature}.
We omit the details and quote only the final result for the Stiefel manifold under the Euclidean metric:
For a tangent plane section spanned by an ONB $\{X=Q\begin{pmatrix}
    A_1\\
    B_2
\end{pmatrix} ,Y= \begin{pmatrix}
    A_2\\
    B_2
\end{pmatrix}\}$,
\begin{eqnarray}
	\nonumber
	K^{St}_e(X,Y)
	&=& \|B_1A_2 - B_2A_1\|_F^2 
       + \frac12\|B_1B_2^T-B_2B_1^T\|_F^2  - \frac12\|B_1^TB_2 - B_2^TB_1\|_F^2\\  
	\label{eq:seccurv_Stief_euclid}
	& & 
    + \frac14 \| [A_1,A_2] - (B_2^TB_1 - B_1^TB_2)\|_F^2,
\end{eqnarray}
see \cite[Prop 4.2, eq. (34) with $\alpha =1$]{nguyen2022curvature}.
In the same work, it is shown that the interval $[-\frac12, 1]$ is contained in the range of the Euclidean curvature and it is conjectured that $-\frac12\leq K^{St}_e(X,Y) \leq 1$.
We confirm this conjecture.
The Euclidean sectional curvature formula shares many terms with the canonical one from \eqref{eq:seccurv_Stief} but with different factors. The occurrence of a negative contribution is annoying in view of estimates.
As before, we restrict our considerations to plane sections spanned by ONBs. Shifting the base point has no impact on the curvature. Hence, w.l.o.g. $Q=I$, $\{X=\begin{pmatrix}
    A_1\\
    B_2
\end{pmatrix} ,Y= \begin{pmatrix}
    A_2\\
    B_2
\end{pmatrix}\} \subset T_{[I]}St(n,p)$.
The matrices $X,Y$ form a (Riemannian) ONB if and only if
\[
 1 = \|A_i\|_F^2 + \|B_i\|_F^2, \hspace{0.2cm} i=1,2, \quad 0 = \tr(A_1^TA_2) + \tr(B_1^TB_2).
\]

\paragraph{The lower bound on the Euclidean Stiefel curvature}
For an ONB $X,Y$ consisting of subblocks $A_i,B_i$ as above, let $\alpha_i = \|A_i\|_F$, $\beta_i= \|B_i\|_F$, $i=1,2$ so that $\alpha_i^2+ \beta_i^2 = 1$. 
By the reversed triangle inequality,
\begin{eqnarray*}
    K^{St}_e(X,Y) &\geq&
     \|B_1A_2 - B_2A_1\|_F^2 
       + \frac12\|B_1B_2^T- B_2B_1^T\|_F^2  - \frac12\|B_1^TB_2 - B_2^TB_1\|_F^2\\  
    &&+ \frac14 \biggl(\| [A_1,A_2]\|_F^2 +  \|B_2^TB_1 - B_1^TB_2\|_F^2\\
    && -2\| [A_1,A_2]\|_F\|B_2^TB_1 - B_1^TB_2\|_F\biggr)\\
    &=&  \|B_1A_2 - B_2A_1\|_F^2 + \frac14 \| [A_1,A_2]\|_F^2 - \frac12 \| [A_1,A_2]\|_F\|B_1^TB_2 - B_2^TB_1 \|_F \\
    && + \frac12 \|B_1B_2^T - B_2B_1^T\|_F^2 -\frac14  \| B_1^TB_2 - B_2^TB_1\|_F^2\\
    &\geq& -\frac12 \| [A_1,A_2]\|_F\|B_1^TB_2 - B_2^TB_1 \|_F
           -\frac14  \|B_1^TB_2- B_2^TB_1 \|_F^2\\
    &\geq& -\frac{1}{2} \alpha_1\alpha_2\sqrt{2}\beta_1\beta_2 - \frac12 \beta_1^2\beta_2^2\\
    &=& -\frac12 \left(\beta_1\beta_2\left(\sqrt{2}\sqrt{1-\beta_1^2}\sqrt{1-\beta_2^2} + \beta_1\beta_2\right)\right).
\end{eqnarray*}
where we have used the Wu-Chen matrix inequality \eqref{eq:WuChen_revisited}.
The minimum of this function in the admissible range of $(\beta_1,\beta_2)\in [0,1]^2$ is attained at the corner point $\beta_1=1=\beta_2$, so that
\[
  K^{St}_e(X,Y)  \geq -\frac12.
\]
 As in the case of the canonical metric, the extrema are attained for tangent vectors with zero $A$-blocks.
For a section spanned by such an ONB, the curvature formula reduces to
\begin{equation}
\label{eq:reduced_eucl}
K^{St}_e(X,Y) =
     \frac12\|B_1B_2^T-B_2B_1^T\|_F^2  - \frac14\|B_1^TB_2 - B_2^TB_1\|_F^2.  
\end{equation}
The lower bound is attained for a normalized, orthogonal matrix pair $B_1,B_2$ such that the first positive term vanishes, while the second negative term becomes extremal. A corresponding matrix pair can be read off from 
 \eqref{eq:max_WuChen14} and \eqref{eq:max_WuChen15}, when choosing $\sigma_1=1=b_*$ and $\sigma_2 = 0 =c_*$.
The associated matrices that span an extremal tangent plane section are
$X= \begin{pmatrix}
    \mathbf{0}\\
    B_1
\end{pmatrix},
Y = \begin{pmatrix}
    \mathbf{0}\\
    B_2
\end{pmatrix}$, where 
\begin{equation}
\label{eq:St_eucl_curvmin_matrices}
    B_1 = \begin{pmatrix}
    0&1\\
    0&0
    \end{pmatrix},
    \quad 
    B_2 = \begin{pmatrix}
    1&0\\
    0&0
    \end{pmatrix}.
\end{equation}
For these matrices, $ \frac12\|B_1^TB_2 - B_2^TB_1\|_F^2 = 1$, while $ \frac12\|B_1B_2^T - B_2B_1^T\|_F^2 = 0$.
Note that this matrix pair is a maximizer of the trace term in \Cref{lem:trace_max}, item {\em (1)}.
%
%
\paragraph{The upper bound on the Euclidean Stiefel curvature}
Continuing the above notation and introducing $\gamma = \|B_1^TB_2 - B_2^TB_1\|_F\in [0, \sqrt{2}\beta_1\beta_2]$, we obtain for the upper bound on the Euclidean Stiefel curvature
\begin{eqnarray}
\nonumber
    K^{St}_e(X,Y) &\leq&
    \|B_1A_2 - B_2A_1\|_F^2 
       + \frac12\|B_1B_2^T- B_2B_1^T\|_F^2  - \frac12\|B_1^TB_2 - B_2^TB_1\|_F^2\\ 
\nonumber
    &&+ \frac14 \biggl(\| [A_1,A_2]\|_F^2 +  \| B_1^TB_2-B_2^TB_1\|_F^2\\
\nonumber
    && \hspace{1.0cm}+2\| [A_1,A_2]\|_F\|B_1^TB_2 - B_2^TB_1\|_F\biggr)\\
\nonumber
    &\leq& \frac12 \alpha_2^2\beta_1^2 +\frac12 \alpha_1^2\beta_2^2 + \alpha_1\alpha_2\beta_1\beta_2  + \beta_1^2\beta_2^2 - \frac12 \gamma^2\\
\nonumber
    && +\frac14\left(\alpha_1^2\alpha_2^2 + \gamma^2 + 2\alpha_1\alpha_2 \gamma \right)\\
\nonumber
    &=& \frac12 (1-\beta_2^2)\beta_1^2 +\frac12 (1-\beta_1^2)\beta_2^2 + \sqrt{1-\beta_1^2}\sqrt{1-\beta_2^2}\beta_1\beta_2  + \beta_1^2\beta_2^2 \\
\label{eq:euclid_curv_est}
    && - \frac12 \gamma^2+\frac14\left((1-\beta_1^2)(1-\beta_2^2) + \gamma^2 + 2 \gamma \sqrt{1-\beta_1^2}\sqrt{1-\beta_2^2}\right).
\end{eqnarray}
When considering all other parameters as fixed, the last expression is a parabola in $\gamma$; the parametric maximum is $\gamma_* = \gamma_*(\beta_1,\beta_2) = \sqrt{1-\beta_1^2}\sqrt{1-\beta_2^2}$.
Inserting this $\gamma$ in \eqref{eq:euclid_curv_est} yields
\begin{eqnarray*}
K^{St}_e(X,Y)  &\leq& 
\frac12 (1-\beta_2^2)\beta_1^2 +\frac12 (1-\beta_1^2)\beta_2^2 + \sqrt{1-\beta_1^2}\sqrt{1-\beta_2^2}\beta_1\beta_2  + \beta_1^2\beta_2^2\\
&&+ \frac12 (1-\beta_1^2)(1-\beta_2^2).
\end{eqnarray*}
%
This expression has its global maximum at $\beta_1=1=\beta_2$, which can be verified by elementary means similar to the discussion of \eqref{eq:brute_force_ineq}. We omit the details. Again, the associated matrix maximizers feature a zero $A$-block.
Reconsidering \eqref{eq:reduced_eucl}, we can utilize 
 \eqref{eq:max_WuChen14} and \eqref{eq:max_WuChen15} to find a maximizing matrix pair. This time one needs to  choose $\sigma_1=1=c_*$ and $\sigma_2 = 0 =b_*$.
 (Mind the "transposition duality" that has been mentioned before.)
The associated matrices that span an extremal tangent plane section are
$X= \begin{pmatrix}
    \mathbf{0}\\
    B_1
\end{pmatrix},
Y = \begin{pmatrix}
    \mathbf{0}\\
    B_2
\end{pmatrix}$, where 
\begin{equation}
\label{eq:St_eucl_curvmax_matrices}
    B_1 = \begin{pmatrix}
    0&0\\
    1&0
    \end{pmatrix},
    \quad 
    B_2 = \begin{pmatrix}
    1&0\\
    0&0
    \end{pmatrix}.
\end{equation}
For these matrices, $ \frac12\|B_1^TB_2 - B_2^TB_1\|_F^2 = 0$, while $ \frac12\|B_1B_2^T - B_2B_1^T\|_F^2 = 1$.
Note that this matrix pair is a maximizer of the trace term in \Cref{lem:trace_max}, item {\em (2)}.
Unlike the canonical metric, the Euclidean metric achieves the cases of extreme curvature for tangent plane sections spanned by rank-1 matrices.

The next theorem summarizes the findings, extended by the cases of special dimensions.
\begin{theorem}
	\label{thm:StiefelCurvBound_eucl}
	The sectional curvature under the Euclidean metric on the Stiefel manifold $St(n,p)$, $n\geq p$, is globally bounded by 
	\[
	-\frac 12 \leq K^{St}_e(X,Y) \leq 1.
	\]
 \begin{enumerate}
     \item For $p\geq 2, n\geq p+2$, the bound is sharp.
 Up to trace-preserving transformations, the maximum and minimum curvature is attained for the tangent plane associated with the matrices of \eqref{eq:St_eucl_curvmax_matrices} and
 \eqref{eq:St_eucl_curvmin_matrices}, respectively.
    \item For $p=1$ and $n\geq 3$, the sectional curvature has a constant value of 
    $$K^{St}_e(X,Y)\equiv 1.$$ 
    \item For $p=n\geq 4$, it holds
    	\[
	0 \leq K^{St}_e(X,Y) \leq \frac14.
	\]
       \item For $p=n=3$, it holds
    	\[
	K^{St}_e(X,Y) \equiv \frac18.
	\]
    \item For $(n,p)=(3,2)$ the sectional curvature is bounded by
    \[ 
        -\frac12\leq K^{St}_e(X,Y)\leq \frac12.
    \]
    The bounds are sharp. The lower bound is sharp and is attained for the matrices from \eqref{eq:St_eucl_curvmin_matrices}, but reduced to their first row.
    The upper bound is attained, e.g., for
    \[
    X = \begin{pmatrix}
        0&0\\
        0&0\\
        -1&0
    \end{pmatrix}, \quad 
    Y = \begin{pmatrix}
        0&-\frac{1}{\sqrt{2}}\\
        \frac{1}{\sqrt{2}}&0\\
        0&0
    \end{pmatrix}.
    \]
        \item For $n\geq 4$ and $(n,p)=(n,n-1)$ the sectional curvature is bounded by
    \[ 
        -\frac12\leq K^{St}_e(X,Y)\leq \frac23.
    \]
    The lower bound is sharp and is attained for the matrices from \eqref{eq:St_eucl_curvmin_matrices}, but reduced to their first row.
 \end{enumerate}
 \begin{conjecture}
 For $n\geq 4$ and $(n,p)=(n,n-1)$, sharp bounds on the sectional curvature of the Stiefel manifold under the Euclidean metric are 
    \[ 
        -\frac12\leq K^{St}_e(X,Y)\leq \frac12.
    \]
 \end{conjecture}

\end{theorem}
\begin{proof}(of \Cref{thm:StiefelCurvBound_eucl})
The global bounds have already been established by considerations that are independent of the dimensions. However, the bounds are not sharp in all dimensions.\\
    On 1.: In the cases, $p\geq 2, n\geq p+2$, both extremal pairs from 
    \eqref{eq:St_eucl_curvmax_matrices} and
 \eqref{eq:St_eucl_curvmin_matrices} fit into the tangent space and thus provide examples in which the global lower and upper bounds are reached.\\
   On 2.: For $p=1$ and $n\geq 3$, the $A$-block in a tangent vector is a skew-symmetric $(1\times 1)$--matrix, i.e., $A=\mathbf{0}$. Moreover,
   the $B$-blocks are $((n-1)\times 1)$--column vectors. Hence, $\|B_1^TB_2 - B_2^TB_1\|_F=0$.
   The curvature formula reduces to
   \[
    	K^{St}_e(X,Y) = \frac12\|B_1B_2^T-B_2B_1^T\|_F^2.
   \]
   Because $X,Y$ form an ONB, we have in this case $B_i^TB_i = 1$ and $B_1^TB_2 = 0$. As a consequence, $\frac12\|B_1B_2^T - B_2B_1^T\|_F^2 = 
   \frac12(2\|B_1\|^2_2\|B_2\|_2^2 - 2\langle B_1,B_2 \rangle_2^2) = 1$.\\
   On 3.: For $p=n\geq 4$, the $B$--blocks of the spanning tangent vectors vanish.
   The curvature formula becomes
   \[
    	K^{St}_e(X,Y) = \frac14\|[A_1, A_2]\|_F^2.
   \]
   Now, the upper bound of $\frac14$ is from \eqref{eq:comm_skew4},  while the lower bound of $0$ is obvious.
   The upper bound is attained, e.g., for the $(4\times 4)$ skew-symmetric matrices $A_1,A_2$ listed at the beginning of \Cref{sec:numex_3}. For a full discussion of the sharpness of the bound, see \cite{Ge2014}.
   The lower bound is attained for any orthogonal, commuting pair of skew-symmetric matrices,
   e.g.,
   \[
   A_1 = \frac{1}{\sqrt{2}} 
   \begin{pmatrix}[cc|cc]
   0 & 1 & 0 & 0\\
   -1 &0&0&0\\\hline
      0&0&0&0\\
   0&0&0&0
   \end{pmatrix},\quad 
   A_2 = \frac{1}{\sqrt{2}}
   \begin{pmatrix}
     [cc|cc]
   0 & 0 & 0 & 0\\
   0 &0&0&0\\\hline
   0&0&0&1\\
   0&0&-1&0  
   \end{pmatrix}.
   \]
   Of course, the curvature is the same for all matrices of higher dimensions that feature the above examples as sub-blocks and are otherwise filled up with zeros.\\
   On 4.: This is in analogy to item 4. of \Cref{thm:StiefelCurvBound}. 
   The argument is the same as outlined in \Cref{app:curv_special_dim}.\\
   On 5.: For $(n,p) = (3,2)$, the matrix blocks $A_1, A_2\in\Skew(2)$ and $B_1, B_2\in\R^{1\times 2}$ can be written as 
   \[A_1 = \alpha\begin{pmatrix}[cc] 0&-1\\1&0\end{pmatrix},\;A_2 = a\begin{pmatrix}[cc] 0&-1\\1&0\end{pmatrix},\;B_1 = \begin{pmatrix}[cc] \beta_1 & \beta_2\end{pmatrix},\;B_2 = \begin{pmatrix}[cc] b_1 & b_2\end{pmatrix},\]
   for $\alpha,\beta_1,\beta_2,a,b_1,b_2\in\R$. 
   The curvature formula reduces to
   \begin{eqnarray}
	\nonumber
	K^{St}_e(X,Y)
	&=& \|B_1A_2 - B_2A_1\|_F^2 
        - \frac14\|B_1^TB_2 - B_2^TB_1\|_F^2
    \end{eqnarray}
    and the orthogonality constraints 
    \[
 1 = \|A_i\|_F^2 + \|B_i\|_F^2, \hspace{0.2cm} i=1,2, \quad 0 = \tr(A_1^TA_2) + \tr(B_1^TB_2)
\] 
    translate to
    \begin{eqnarray}
	\nonumber
    \alpha^2 = \frac12 - \frac12(\beta_1^2 + \beta_2^2), \quad a^2 = \frac12 - \frac12(b_1^2 + b_2^2), \quad \beta_1 b_1 + \beta_2 b_2 = -2\alpha a.
    \end{eqnarray}
    By squaring both sides of the equation $\beta_1 b_1 + \beta_2 b_2 = -2\alpha a$ and exploiting the equations for $\alpha^2$ and $a^2$, we obtain
    \begin{eqnarray}
	\nonumber
	(\beta_1^2 + \beta_2^2) + (b_1^2 + b_2^2) - 1 = (\beta_1 b_2 - \beta_2 b_1)^2 \geq 0
    \end{eqnarray}
    and therefore $(\beta_1^2 + \beta_2^2) + (b_1^2 + b_2^2) \geq 1$. By exploiting the equations resulting from the orthogonality constraints once again, the bound on the curvature is obtained:
    \begin{eqnarray}
	\nonumber
	K^{St}_e(X,Y)
	&=& \frac32 - \big((\beta_1^2 + \beta_2^2) + (b_1^2 + b_2^2)\big) \leq \frac12.
    \end{eqnarray}
    On the other hand, $-\frac12\leq K^{St}_e(X,Y)$ applies because $(\beta_1^2 + \beta_2^2) + (b_1^2 + b_2^2)\leq 2$.
\\
For an alternative, but more involved argument, see \cite[Prop. 6.1]{nguyen2022curvature}, which gives tight bounds in the case $p=2$ for a family of metrics including the Euclidean and the canonical one.\\
   On 6.: For $p=n-1$, the matrix blocks $B_1, B_2\in \R^{(n-(n-1))\times (n-1)}$ are row vectors in $\R^{1\times (n-1)}$.
   The global lower curvature bound of $K_e\geq -\frac12$ continues to hold and the sharpness of the bound is confirmed by the example of the tangent space section associated with the matrix blocks from \eqref{eq:St_eucl_curvmin_matrices}, which fit into this setting dimension-wise, when reduced to their first row and possible filled up with zeros.
By taking 
   $\|B_1B_2^T - B_2B_1^T\|_F^2 = 0$ into account,
   the curvature formula reads
\begin{eqnarray}
	\nonumber
	K^{St}_e(X,Y)
	&=& \|B_1A_2 - B_2A_1\|_F^2 
        - \frac12\|B_1^TB_2 - B_2^TB_1\|_F^2\\  
	& & \nonumber
    + \frac14 \| [A_1,A_2] - (B_2^TB_1 - B_1^TB_2)\|_F^2.
\end{eqnarray}
   The upper bound of $K_e\leq\frac23$ can be established in the same manner as before.
   We omit the details, on the one hand because this is straightforward, on the other hand, because we do not believe that this bound is sharp anyways.
\end{proof}
\subsection{Impact on the injectivity radius of the Stiefel manifold}
\label{sec:injectivity_radius}
On a Riemannian manifold $\mathcal{M}$, the injectivity radius is the largest possible radius within which geodesics are unique and lengths-minimizing, regardless of where you start from.
In loose words, as long as you stay within the injectivity radius when travelling along a geodesic, you are guaranteed not to travel unnecessary distances.
This concept is formalized with introducing the Riemannian exponential map at a point $p\in\mathcal{M}$, which sends a tangent vector $v\in T_p\mathcal{M}$, to the endpoint of the geodesic on the unit interval that starts from $p$ with velocity $v$.
For a precise definition of the notion of the injectivity radius, we refer the reader to \cite[Chap. 13]{DoCarmo2013riemannian}.

A classical result from Riemannian geometry relates the injectivity radius with the sectional curvature:
\begin{theorem}[Klingenberg, stated as Lemma 6.4.7 in \cite{petersen2016riemannian}]
\label{thm:Klingenberg}
	Let $\mathcal{M}$ be a compact Riemannian manifold with sectional curvatures bounded by $K(X,Y)\leq C$, where $C>0$. Then the injectivity radius $\mathrm{inj}(p)$ at any $p\in \mathcal{M}$ satisfies
	\[
	\mathrm{inj}(p)\geq \min\left\{\frac{\pi}{\sqrt{C}}, \frac12 l_p\right\},
	\]
	where $l_p$ is the length of a shortest closed geodesic starting from $p$.
	For the global injectivity radius, it holds
		\[
	\mathrm{inj}(\mathcal{M})  \geq\frac{\pi}{\sqrt{C}}
	\quad \text{ or } \quad \mathrm{inj}(\mathcal{M}) =\frac12 l,
	\]
	where $l$ is the length of a shortest closed geodesic on $\mathcal{M}$.
\end{theorem}
From \Cref{thm:StiefelCurvBound} and \Cref{thm:StiefelCurvBound_eucl}, one obtains for $p\geq 2, n\geq p+2$.
\begin{enumerate}
\item $\mathrm{inj}(St(n,p))\geq \sqrt{\frac{4}{5}}\pi$ (canonical metric).
\item $\mathrm{inj}(St(n,p)) = \pi$ unless there is a closed geodesic of length strictly smaller than $2\pi$ (Euclidean metric).
\end{enumerate}
Under the canonical metric, it is clear that closed geodesics have a length of at least $2\pi$, see \cite[Chapter 5]{Rentmeesters2013} or the recent preprints \cite{absil2024ultimate}, \cite{stoye2024injectivity}.
A closed geodesic on the unit interval $[0,1]$ of length $2\pi$ on the Stiefel manifold $St(4,2)$ under the Euclidean metric is
\[
    c(t) = \exp_m\left(t\begin{pmatrix}
        2A & -B^T\\
        B  & \mathbf{0} 
    \end{pmatrix}\right)
    \begin{pmatrix}
      \exp_m(-tA)\\
      \mathbf{0}
    \end{pmatrix}, \quad A = \mathbf{0}\in \Skew(2), 
    \quad B = \begin{pmatrix}
        2\pi & 0\\
        0 & 0
    \end{pmatrix}.
\]
The form of this geodesic follows the representation from \cite[Prop. 1]{ZimmermannHueper2022}.
Note that $c(0) = c(1) = 
\begin{pmatrix}
I_2\\
\mathbf{0}
\end{pmatrix}$ and that the length of $c|_{[0,1]}$ is $l=2\pi$.
Of course, this geodesic can be embedded in higher dimensions of $p\geq2, n\geq p+2$.
Examples of closed geodesics and their length with respect to the metrics of the family from 
\cite{HueperMarkinaLeite2020} are given in \cite[Section 6]{absil2024ultimate}.
We believe that the above example is a shortest closed geodesic under the Euclidean metric.
A more detailed discussion is beyond the scope of this paper, but will be provided in \cite{stoye2024injectivity}.
Closely related investigations can be found in \cite{absil2024ultimate}.
\section{Numerical Experiments}
\label{sec:numex}
In this section, we illustrate the behavior of the sectional curvature on $SO(n), Gr(n,p)$ and $St(n,p)$ at special parametric sections, where the rank of the spanning matrices increase with the parameter. 
We also investigate the behavior of the generic sectional curvature on these manifolds when increasing the dimension $p$. For $St(n,p)$ both the canonical and the Euclidean metric feature in the experiments.
\subsection{Experiment 1}
\label{sec:Exp1}
\begin{figure}[ht]
	\centering
	\includegraphics[width=0.9\textwidth]{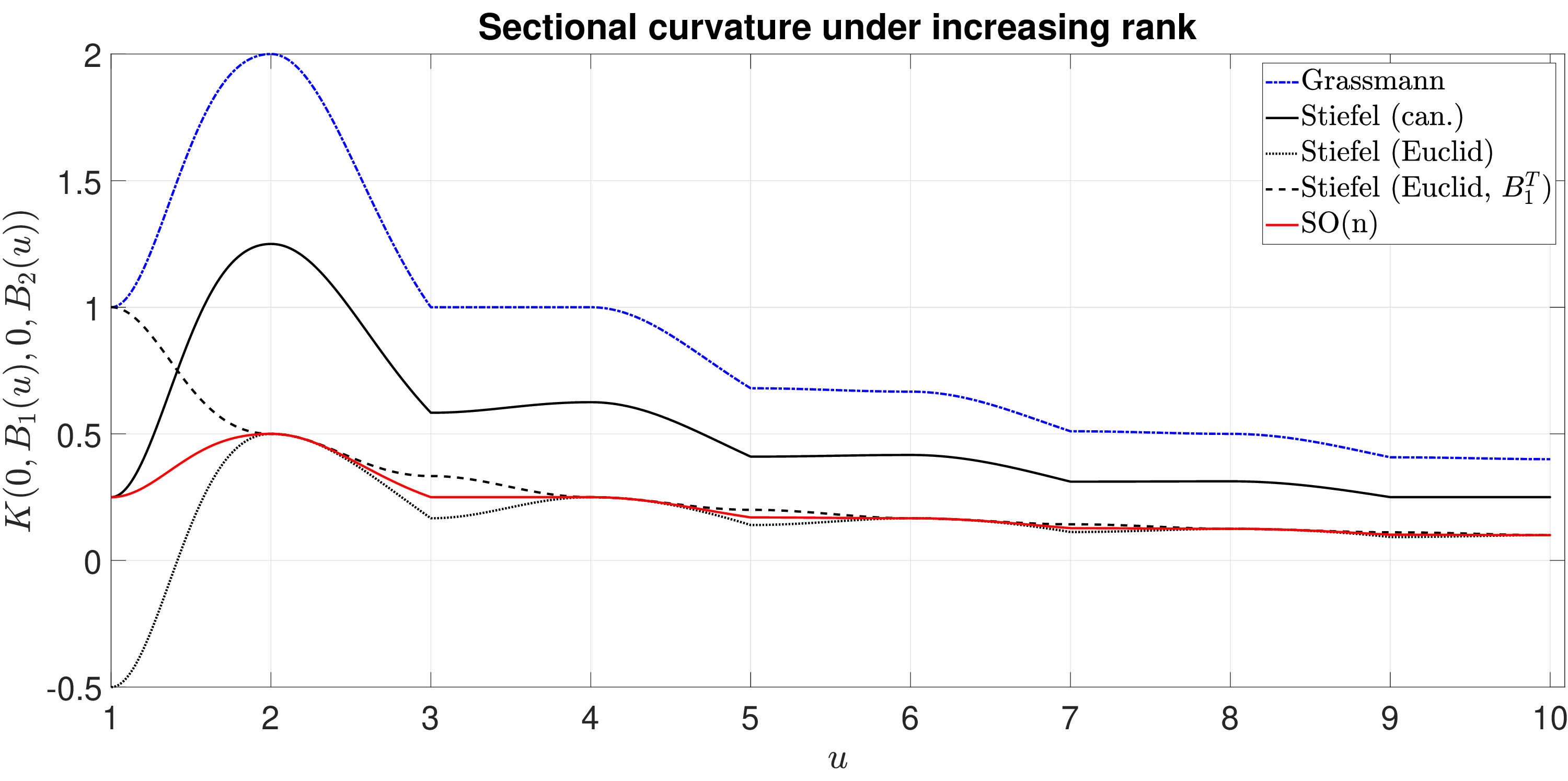}
	\caption{(Corresponding to \Cref{sec:Exp1}.) Sectional curvature on $SO(20)$, $St(20,10)$, $Gr(20,10)$ for the tangent sections defined by \eqref{eq:exp1_matrices}. For $SO(20)$, $Gr(20,10)$, and $St(20,10)$ under the canonical metric, the respective global sectional curvature maximum is attained for the matrices with subblocks $B_1$, $B_2$ as stated in \Cref{thm:StiefelCurvBound}.
   For $St(n,p)$ under the Euclidean metric, the extreme curvature cases occur for the rank-one matrices from \eqref{eq:St_eucl_curvmax_matrices} (max) and \eqref{eq:St_eucl_curvmin_matrices} (min).
	}
	\label{fig:exp1}
\end{figure}

In the first experiment, we start with the rank-one matrices
\[
  B_1 = \begin{pmatrix}[c|c]
	\begin{matrix} 0 & 1\\ 0 & 0 \end{matrix} & \mathbf{0} \\
	\hline
	\mathbf{0}& \mathbf{0} \\ 
\end{pmatrix}
\mbox{ and }
B_2 = 
\begin{pmatrix}[c|c]
	\begin{matrix} 1  & 0\\ 0 & 0 \end{matrix} & \mathbf{0} \\
	\hline
	\mathbf{0}& \mathbf{0} \\ 
\end{pmatrix} \in \R^{10\times 10}.	
\]
Then we fill the remaining diagonal entries of $B_2$ and the super- and sub-diagonal entries of $B_1$ one after the other
in the following way
\[
  B_1(u) = \begin{pmatrix}[c|c|c]
		\begin{matrix} 0 & 1\\ -u_2 & 0 \end{matrix} & &\\ \hline
		                                          & \ddots &\\ \hline
		                                          &        & \begin{matrix} 0 & u_9\\ -u_{10} & 0 \end{matrix}
\end{pmatrix}
\mbox{ and }
B_2(u) = 
\begin{pmatrix}[cccc]
                   1  &       &        & \\
                      & u_2   &        &\\
	                  &       & \ddots & \\
	                  &       &        & u_{10} 
	\end{pmatrix} \in \R^{10\times 10}.
\]
In the beginning, $u_2=u_3=\ldots=u_{10} = 0$.
Then, we first increase $u_2\in[0,1]$ linearly until the upper bound $1$ is reached.
Then, we keep $u_2=1$ and let $u_3$ run through $[0,1]$ and so on.
For each matrix pair $B_1(u),B_2(u)$ under this procedure, we compute the tangent vectors
\begin{equation}
	\label{eq:exp1_matrices}
	X(u) = \begin{pmatrix}
		0 & -B_1(u)^T\\
	B_1(u)& 0
	\end{pmatrix}, \quad
		Y(u) = \begin{pmatrix}
		0 & -B_2(u)^T\\
		B_2(u)& 0
	\end{pmatrix}
\end{equation}
and the sectional curvatures
associated with the tangent planes spanned by $X(u),Y(u)$ for the special orthogonal group, the Stiefel manifold and the Grassmann manifold,
\[
	K^{SO}(X(u),Y(u)),\quad K^{St}_c(X(u),Y(u)), \quad K^{St}_e(X(u),Y(u)),\quad K^{Gr}(X(u),Y(u)).
\]
It is understood that the matrices $X(u),Y(u)$ are normalized according to the metric of the respective manifold. To catch both the maximal and the minimal case of the Stiefel curvature under the Euclidean metric, the experiment is performed twice. Once with $B_1(u)$ as above and once with $B_1^T(u)$. 
The results are displayed in \Cref{fig:exp1}.
The figure confirms that for the matrices under consideration, the global curvature maximum on all manifolds under consideration is reached for tangent planes spanned by $X,Y$ that feature zero skew-symmetric $A$-blocks.
For $SO(n), Gr(n,p)$ and $St(n,p)$ under the canonical metric, 
the maximum sectional curvature occurs for the rank-two subblocks $B_1,B_2$ from \eqref{eq:secmaximizers_n4p2}.
Under the Euclidean metric the curvature maximum and minimum occur for the rank-one matrices $B_1,B_2$ from \eqref{eq:St_eucl_curvmax_matrices} and \eqref{eq:St_eucl_curvmin_matrices}, respectively.
\subsection{Experiment 2: Curvature associated with pseudo-random sections}
\label{sec:Exp2}
\begin{figure}[ht]
	\centering
	
	\includegraphics[width=1.0\textwidth]{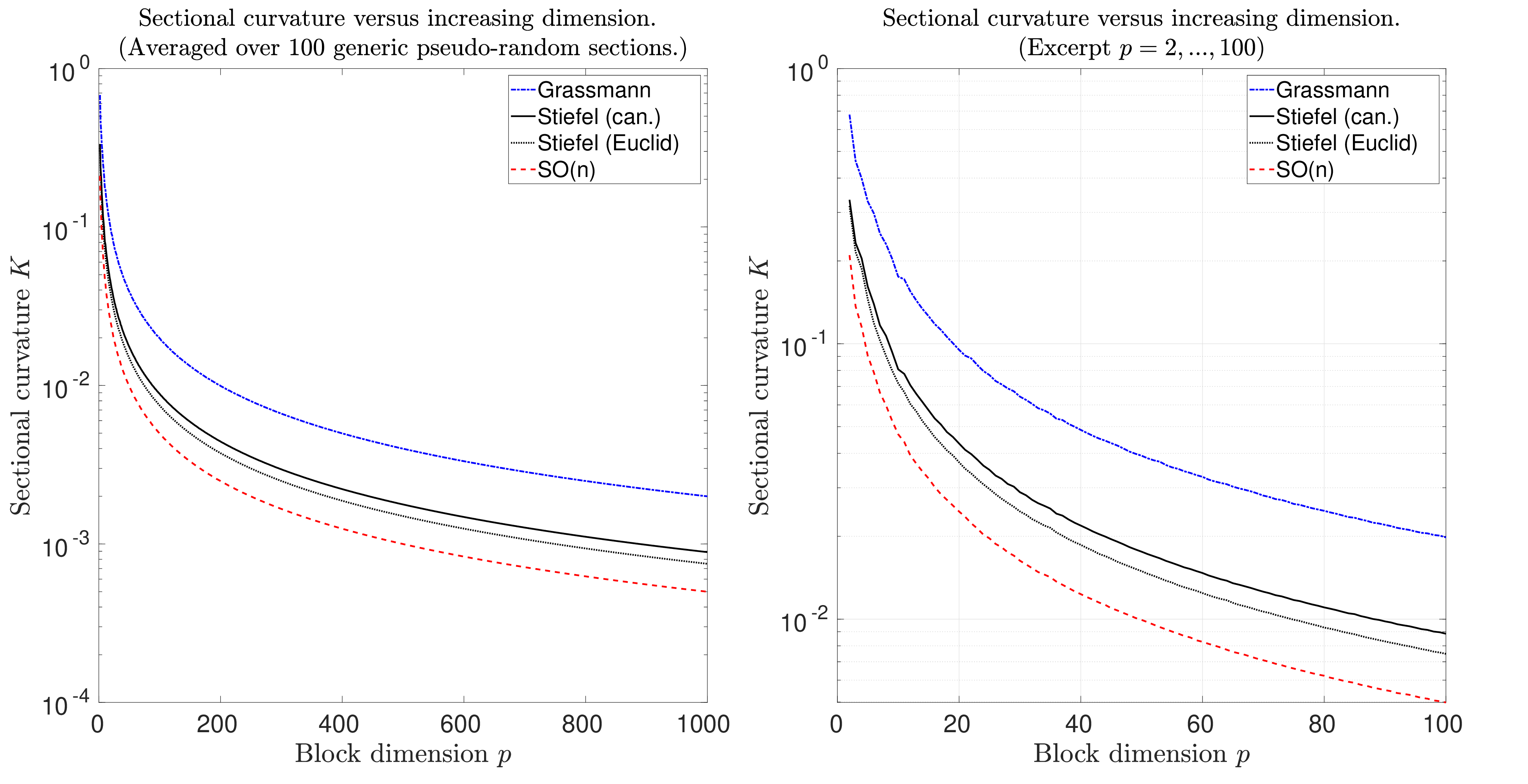}
	\caption{(Corresponding to \Cref{sec:Exp2}.) Averaged sectional curvature for random tangent sections
    $X=X(A_1,B_1,C_1), Y=Y(A_2,B_2,C_2)\in \Skew(2p)$.
    The dimension of the sub-blocks is $p$.
    For Grassmann and Stiefel (canonical metric), $X,Y$ are projected onto the respective horizontal space. For Stiefel under the Euclidean metric, 
    $X= \begin{pmatrix}
        A_1\\
        B_1
    \end{pmatrix}, 
    Y=\begin{pmatrix}
        A_2\\
        B_2
    \end{pmatrix}
    $
    are formed. In all cases, the tangent vectors are orthonormalized according to the respective metric.
	}
	\label{fig:exp2}
\end{figure}

Next, we run an experiment with random-matrices of increasing dimension.
To this end, we create random skew-symmetric matrices 
\[
	X(A_1,B_1,C_1) = \begin{pmatrix}
	A_1 & -B_1^T\\
	B_1 & C_1
\end{pmatrix},\quad 
    Y(A_2,B_2,C_2)= \begin{pmatrix}
	A_2 & -B_2^T\\
	B_2 & C_2
\end{pmatrix}	\in \Skew(2p).
\]
We start from $p=2$. For each $p=2,\ldots, 1000$, a number of $100$ random pairs $X,Y\in\Skew(2p)$ are computed
and the sectional curvature of 
\begin{itemize}
\item the tangent plane spanned by $X,Y\in T_ISO(2p)$, 
\item the tangent plane spanned by $X_{\mathfrak{m}^{St}}= X(A_1,B_1,\mathbf{0}),
Y_{\mathfrak{m}^{St}}= Y(A_2,B_2,\mathbf{0})
\in T_{[I]}St(2p,p)$ (canonical metric),
\item the tangent plane spanned by $X_{\text{eucl}}=(A_1^T,B_1^T)^T, Y_{\text{eucl}}=(A_2^T,B_2^T)^T\in T_{ (I,0)^T}St(2p,p)$ (Euclidean metric),
\item and the tangent plane spanned by 
$X_{\mathfrak{m}^{Gr}}= X(\mathbf{0},B_1,\mathbf{0})$,
$Y_{\mathfrak{m}^{Gr}}= Y(\mathbf{0},B_2,\mathbf{0}),
\in T_{[I]}Gr(2p,p)$,
\end{itemize}
respectively, is computed.
The result is then averaged over the number of $100$ random runs. 
\Cref{fig:exp2} displays the sectional curvature versus the block dimension $p$ on a logarithmic scale.
It can be seen that the curvature of these "random" planes is largest for $p=2$ and drops considerably, when $p$ is increased. To be precise, for the canonical Stiefel data, the averaged random curvature at $p=2$ is $K=0.33$, while it is $K=8.8\cdot 10^{-3}$ for $p=100$ and  $K=8.9\cdot 10^{-4}$ for $p=1000$.

\subsection{Experiment 3: Impact of the blocks A and B}
\label{sec:numex_3}

\begin{figure}[ht]
		\centering
		\includegraphics[width=1.0\textwidth]{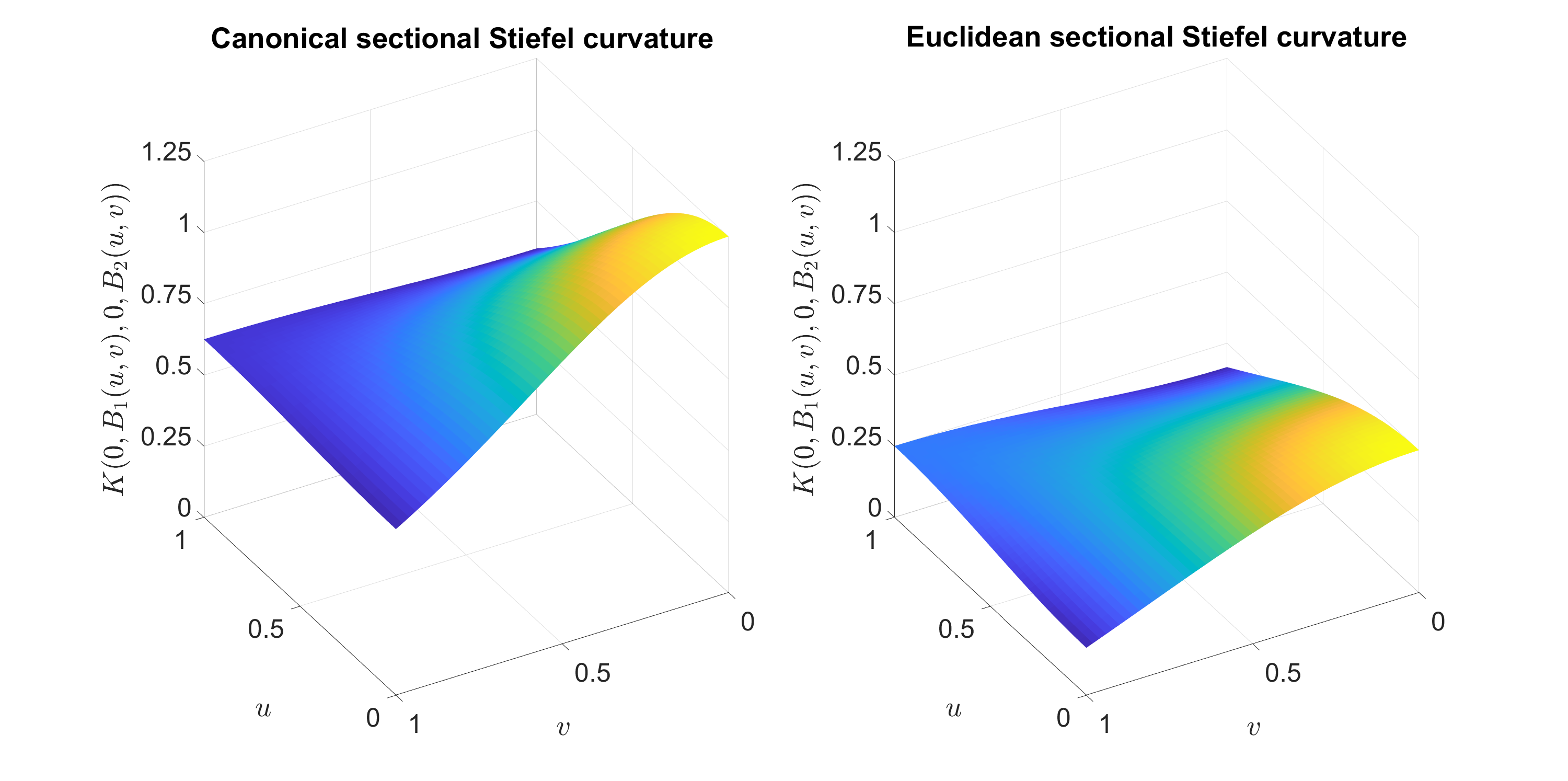}
		\caption{Sectional curvature on Stiefel for the tangent sections spanned by the matrix blocks of \Cref{sec:numex_3}.
     Left: canonical metric, Right: Euclidean metric.
		}
		\label{fig:exp3a}
	\end{figure}

\begin{figure}[ht]
		\centering
		\includegraphics[width=1.0\textwidth]{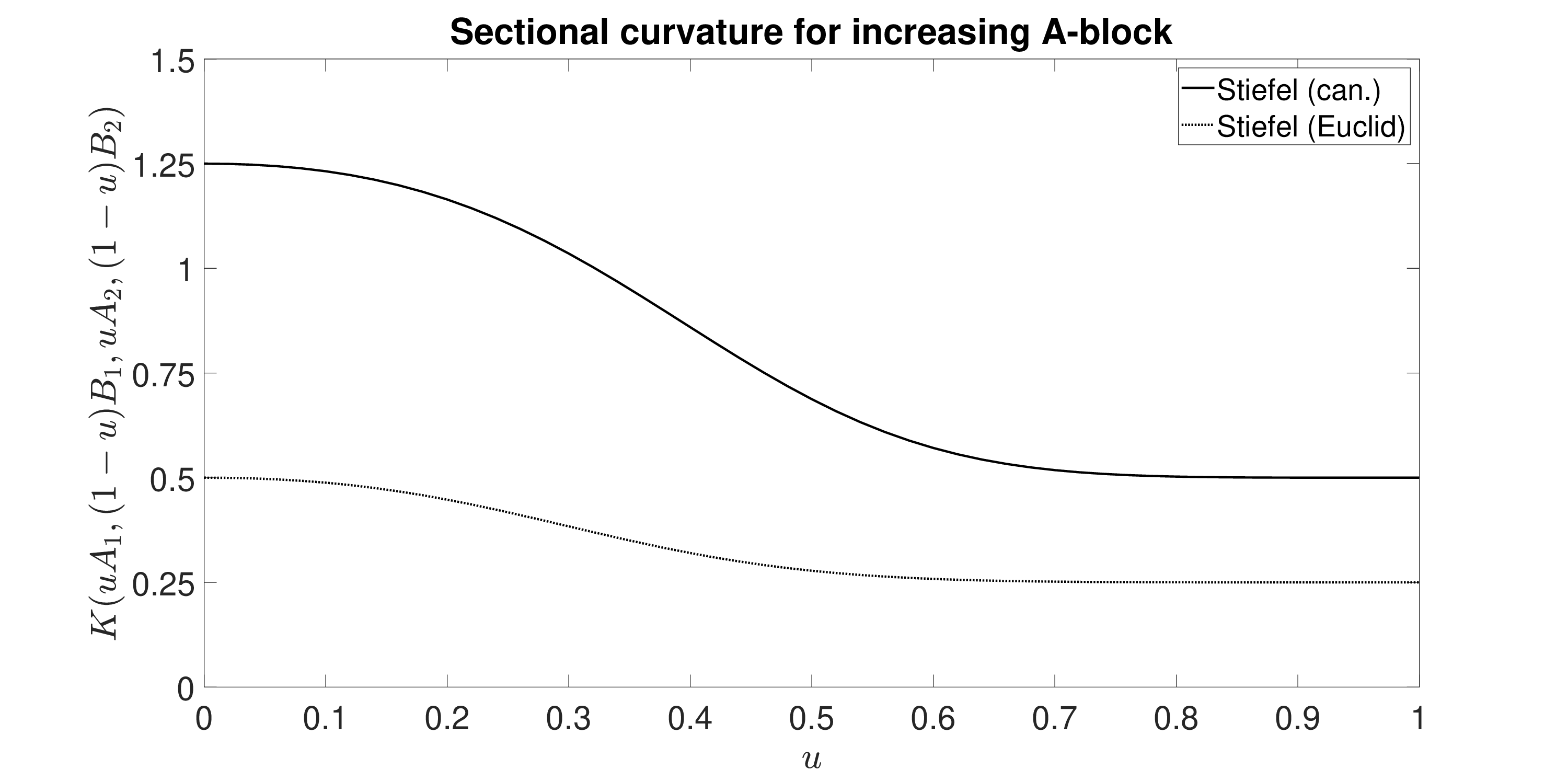}
		\caption{Sectional curvature on Stiefel when the weight in the spanning normalized tangent matrices is shifted from the $B$-blocks to the $A$-blocks.
		}
		\label{fig:exp3b}
	\end{figure}

	Consider the special matrices
	\begin{align*}
		A_1 &= 
		\begin{pmatrix}[cccc]
			0 & 1       & 0       & 0\\
			-1& 0   & 0       &0\\
			0&  0     & 0 &1 \\
			0& 0      & -1        & 0 
		\end{pmatrix},
		\quad	
		&A_2 = 
		\begin{pmatrix}[cccc]
			0 & 0       & 0       & 1\\
			0& 0   & 1       &0\\
			0&  -1     & 0 &0 \\
			-1& 0      & 0        & 0 
		\end{pmatrix}
		\\
		B_1(u,v) &= 
		\begin{pmatrix}[cccc]
			0 & 1       & 0       & 0\\
			-1& 0   & 0       &0\\
			0&  0     & 0 &u \\
			0& 0      & -v        & 0 
		\end{pmatrix},
		\quad	&B_2(u,v) =	
		\begin{pmatrix}[cccc]
			1 & 0       & 0       & 0\\
			0& 1   & 0       &0\\
			0&  0     & u &0 \\
			0& 0      & 0        & v 
		\end{pmatrix}
	\end{align*}
	The matrix pair $A_1,A_2$ maximizes the commutator norm $\|[A_1,A_2]\|_F$.
	For $(u,v)=(0,0)$, the matrix pair $B_1,B_2$ makes the Wu-Chen inequality sharp.
    Let
    \begin{align*}
        K^{St}_c(A_1,B_1,A_2,B_2)&:= 
        K^{St}_c(\begin{pmatrix}
            A_1 &-B_1^T\\
            B_1 &  \mathbf{0}
        \end{pmatrix},
        \begin{pmatrix}
            A_2 &-B_2^T\\
            B_2 &  \mathbf{0}
        \end{pmatrix}),\\
        K^{St}_e(A_1,B_1,A_2,B_2)&:= 
        K^{St}_e(\begin{pmatrix}
            A_1 \\
            B_1
        \end{pmatrix},
        \begin{pmatrix}
            A_2\\
            B_2
        \end{pmatrix}).
    \end{align*}
	Figure \ref{fig:exp3a} displays the function
	$$[0,1]^2\to \R, (u,v) \mapsto K^{St}_{m}(\mathbf{0},B_1(u,v), \mathbf{0} , B_2(u,v)), \quad m\in\{c,e\}$$
 under the canonical metric and the Euclidean metric.
 In both cases, the curvature decreases, when the lower subblocks $B_1, B_2$ get filled up.
	Figure \ref{fig:exp3b} displays the function
	$$[0,1]\to \R, u \mapsto K^{St}_{m}(uA_1,(1-u)B_1(0,0), uA_2 , (1-u)B_2(0,0)), \quad m\in\{c,e\}$$
 under the canonical metric and the Euclidean metric.
 The figure illustrates that the sectional curvature decreases in both cases, when the weight in the normalized tangent vectors is shifted from the $B$-blocks to the $A$-blocks.

%
%
%
%
%
\section{Summary}
\label{sec:conclusions}
This paper resumes the investigation of the sectional curvature on the Stiefel and Grassmann manifolds. In the Grassmann case, tight bounds have been known since the work of Wong \cite{Wong1968} from 1968.
We pay special attention to the maximizers of the curvature bounds and provide refined matrix trace inequalities for this purpose.

An extensive study of the sectional curvature on the Stiefel manifold equipped with a parametric family of Riemannian metrics has been conducted in \cite{nguyen2022curvature}. 
However, since the formulae apply to all members of the family, it can be tedious to extract specific formulae for a particular metric.
Moreover, tight bounds were only obtained for Stiefel manifolds $St(n,p)$ with $p=2$.

Under the canonical metric, which is a member of the parametric family, we confirm prior conjectures regarding the sectional curvature. 
Specifically, we establish that the curvature on the Stiefel manifold equipped with this metric globally does not exceed 5/4. With this addition, we now have a complete account of the curvature bounds in all admissible dimensions.

Under the Euclidean metric, we prove sharp global bounds on the sectional curvature in all admissible dimensions with one single exception:
For the special case $p=n-1$ and $n\geq 4$, we could only establish a sharp lower bound of $K_e\geq-1/2$ and an upper bound of $K_e\leq 2/3$. We believe that the true upper bound is at $K_e\leq 1/2$. We share this believe with the author of \cite{nguyen2022curvature}.

We also show that the sections that maximise the Grassmann curvature are exactly those for which the Stiefel curvature is maximised under the canonical metric, and that these tangent space sections are necessarily spanned by special rank-two matrices.
Under the Euclidean metric, the extreme cases occur for tangent space section spanned by special rank-one matrices.
This supports the observation that 'high curvature means low-rank', which is illustrated by numerical experiments that reveal a decrease in curvature with increasing the rank. 

%
%

%

%
\section*{Acknowledgements}
The authors would like to thank Prof. P.-A. Absil, University of Louvain, for stimulating conversations on the subject.
%
\appendix

\section{The global maximum of the curvature bound of \texorpdfstring{\eqref{eq:brute_force_ineq}}{(22)}}
\label{app:auxresults}
\begin{figure}
    \centering
    \includegraphics[width=0.5\textwidth]{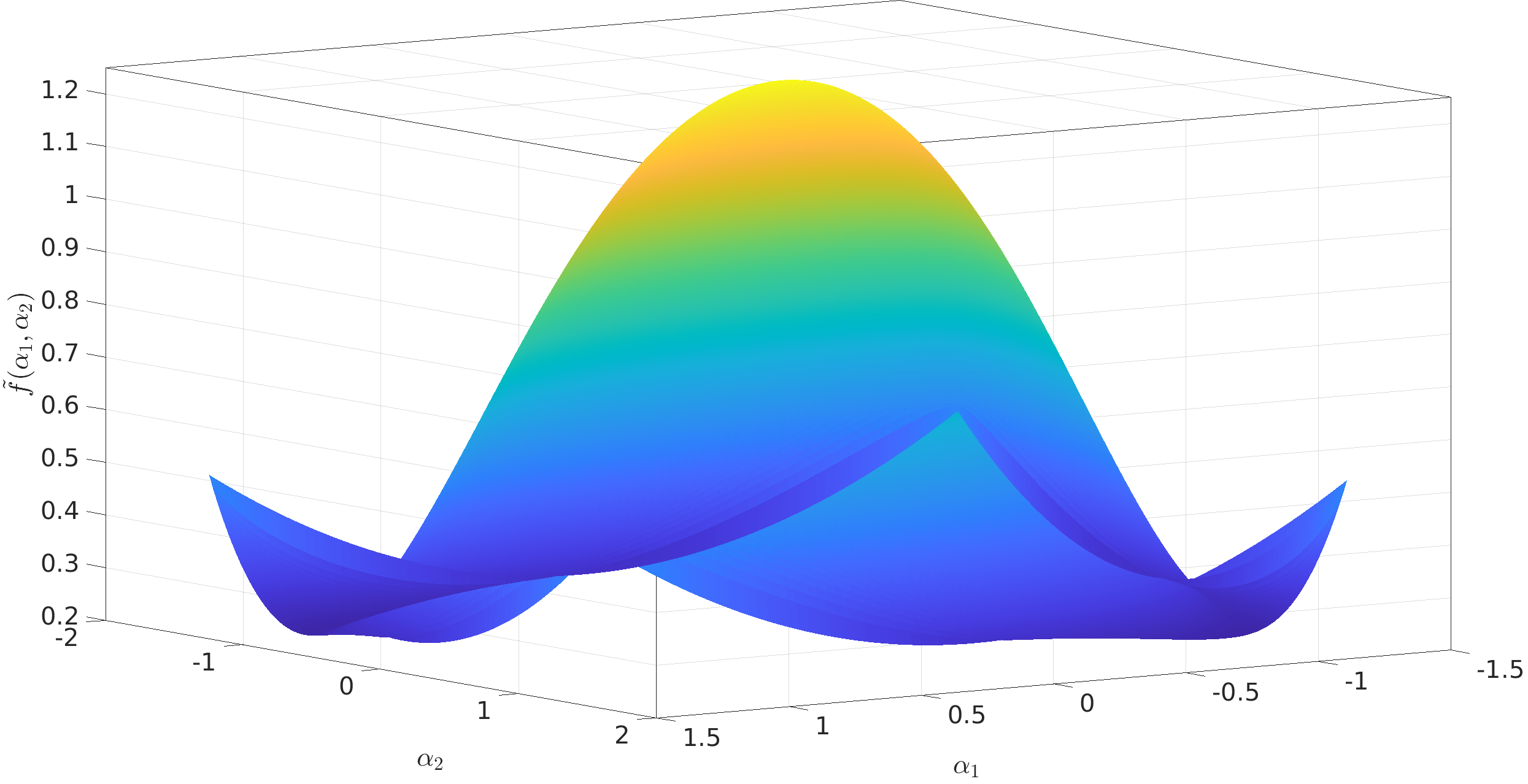}
    \caption{Plot of the function $\tilde f$ from \eqref{eq:globmax}
    that bounds the Stiefel sectional curvature in the range $[-\sqrt{2},\sqrt{2}]^2$.}
    \label{fig:plot_brute_force}
\end{figure}

We formally verify that the function 
in \eqref{eq:brute_force_ineq} that bounds the sectional curvature of the Stiefel manifold 
\begin{equation}
\tilde f: (\alpha_1,\alpha_2)\mapsto \frac54 + \frac{5}{16}\alpha_1^2\alpha_2^2 - \frac12 (\alpha_1^2 + \alpha_2^2)
	+\frac{1+\sqrt{2}}{4}\alpha_1\alpha_2\sqrt{1-\frac12\alpha_1^2}\sqrt{1-\frac12\alpha_2^2}
 \label{eq:globmax}
    \end{equation}
has its global maximum at $(0,0)$ for $(\alpha_1,\alpha_2)\in\left[0,\sqrt{2}\right]^2$.\\
Consider the transformation $\alpha:\left[0,\frac{\pi}{2}\right]\to\left[0,\sqrt{2}\right], r\mapsto 
 \alpha(r) = \sqrt{2}\sin{(r)}$ so that $\sqrt{1-\frac{1}{2}\alpha^2} = \cos{(r)}$ and 
\[
    \alpha\sqrt{1-\frac{1}{2}\alpha^2} = \sqrt{2}\sin{(r)}\cos{(r)} = \frac{\sqrt{2}}{2}\sin{(2r)}.
\]
With parameterizing $\alpha_1=\alpha_1(r),\alpha_2=\alpha_2(s)$ in this form, the task is equivalent to showing that the global maximum of $f(r,s) := \tilde f(\alpha_1(r),\alpha_2(s))$,
\[f(r,s) = \frac54 + \frac{5}{4}\sin(r)^2\sin(s)^2 - \sin(r)^2 - \sin(s)^2 + \frac{1+\sqrt{2}}{8}\sin(2r)\sin(2s),\]
in the admissible range $\left[0,\frac{\pi}{2}\right]^2$ is at $(0,0)$. The gradient of $f$ is
\[\nabla f(r,s) = \begin{bmatrix}
    f_r(r,s)\\
    f_s(r,s)
\end{bmatrix} = \begin{bmatrix}
    \frac{1}{2}\sin(r)\cos(r)\left(5\sin(s)^2 - 4\right) + \frac{1+\sqrt{2}}{4}\cos(2r)\sin(2s)\\
    \frac{1}{2}\sin(s)\cos(s)\left(5\sin(r)^2 - 4\right) + \frac{1+\sqrt{2}}{4}\cos(2s)\sin(2r)
\end{bmatrix}.\]
The condition $f_r(r,s) = 0$ gives
\begin{equation}
    \frac{1}{2}\sin(r)\cos(r)\left(4 - 5\sin(s)^2\right) = \frac{1+\sqrt{2}}{4}\cos(2r)\sin(2s). \label{eq:equalf1v1}
\end{equation}
We assume $r,s\neq\frac{\pi}{4}$ and $\sin(s)^2,\sin(r)^2\neq\frac{4}{5}$ and tackle those special cases afterwards. With excluding those cases we equivalently obtain
\begin{equation}
    \frac{1}{2}\frac{\sin(r)\cos(r)}{\cos(2r)} = \frac{1+\sqrt{2}}{4}\frac{\sin(2s)}{4 - 5\sin(s)^2}. \label{eq:equalf1v2}
\end{equation}
 The left-hand side is greater than zero for $r<\frac{\pi}{4}$ and smaller than zero for $r>\frac{\pi}{4}$ and in each case monotonically increasing. The same holds for the right-hand side with $\sin(s)^2<\frac{4}{5}$ and $\sin(s)^2>\frac{4}{5}$. We only consider the case where both sides are greater than zero. The other cases can be tackled analogously. Due to the monotonicity of both sides, it immediately follows that there is at most one pair $(r,s)$ that fulfills the equation. Suppose $(r,s)$ is a pair satisfying the equation. We will show that in this case $r$ has to be equal to $s$. Assume $r\neq s$. Let $s<r<\frac{\pi}{4}$ and investigate $f_s(r,s) = 0$ which is equivalent to 
\begin{equation}
\frac{1}{2}\frac{\sin(s)\cos(s)}{\cos(2s)} = \frac{1+\sqrt{2}}{4}\frac{\sin(2r)}{4 - 5\sin(r)^2}. \label{eq:equalf2}
\end{equation}
With $s<\frac{\pi}{4}$ we are also in the case where both sides have to be positive in order for the equation to hold. The equation \eqref{eq:equalf2} for $f_s$ is the same as the equation \eqref{eq:equalf1v2} for $f_r$, but with the roles for $r$ and $s$ reversed. Thus, for the pair $(r,s)$ fulfilling the equation \eqref{eq:equalf2}, $r<s$ applies. This contradicts the above assumption. 
The case $r<s$ can be tackled in the same way.
In summary, for the equation to hold, $r$ must be equal to $s$
so that candidates for extrema lie necessarily on the diagonal $r = s$. 
Along this diagonal, the function $\eqref{eq:globmax}$ bounding the sectional curvature becomes a parabola in $\alpha^2$
\begin{equation*}
\tilde f: \alpha\mapsto \frac54 + \frac{3-2\sqrt{2}}{16}\alpha^4 - \frac{3-\sqrt{2}}{4}\alpha^2.
\end{equation*}
It is easy to show that $\tilde f$ is monotonically strictly decreasing for $\alpha\in\left[0,\sqrt{2}\right]$ and therefore the maximum is at $\alpha = 0$.\\
Now, we tackle the remaining cases $r,s=\frac{\pi}{4}$ and $\sin(s)^2,\sin(r)^2=\frac{4}{5}$. At $r=\frac{\pi}{4}$, the condition $f_r(r,s) = 0$ yields $\sin(s)^2 = \frac{4}{5}$. But in this case $f_s(r,s)\neq 0$. The other combinations can be treated analogously (and are also not describing any extreme points).
This completes the analysis and verifies that $\tilde f$ has its unique global maximum in the admissible range $\left[0,\sqrt{2}\right]^2$ at $(0,0)$.
\Cref{fig:plot_brute_force} displays the function $\tilde f$.
%
\section{The sectional curvature of special \texorpdfstring{$St(n,p)$}{St(n,p)}}
\label{app:curv_special_dim}
For the sake of completeness, we give a direct proof of item 4 of \Cref{thm:StiefelCurvBound}:
\begin{quote}
    On $St(3,2)$ and $St(3,3)$, the sectional curvature is $K^{St}(X,Y)\equiv \frac14.$
\end{quote}
The result is already contained in \cite[Prop. 6.1]{nguyen2022curvature}.
It follows also immediately from \Cref{rem:curvSOn_fro} combined with \eqref{eq:comm_skew3}.
\begin{proof}
On $St(3,2)$, tangent vectors are skew-symmetric $(3\times 3)$-matrices,
\[
    X = \begin{pmatrix}
        A & -B^T\\
        B &\mathbf{0}
    \end{pmatrix}
    = \begin{pmatrix}[cc|c]
        0 & -a & -b_1 \\
        a & 0  & -b_2\\
        \hline
        b_1&b_2& 0
    \end{pmatrix}.
\]
For $X\in\Skew(3)$ as above, let  $x:=\text{vec}_3(X)=(a,b_1,b_2)\in\R^3$ be the vector that features the independent entries of $X$ as components.
Let $X,Y\in \Skew(3)$ and $x=\text{vec}_3(X)$, $y=\text{vec}_3(Y)$.
The following facts are well-known and may be verified by elementary means:
\begin{itemize}
    \item $\text{vec}_3([X,Y])= x\times y$.
    \item $\langle X,Y\rangle_F= 2\langle x,y\rangle_2$.
    In particular,  $\|X\|_F^2 = 2\|x\|_2^2$ and $(X\bot_F Y) \Leftrightarrow (x\bot_F y)$.
\end{itemize}
Let $X,Y\in \Skew(3)$ be an ONB w.r.t. the Frobenius inner product.
In this case, $\|x\|_2=\frac{1}{\sqrt{2}}=\|y\|_2$ for the associated vector representations.
The curvature formula of \Cref {rem:curvSOn_fro} applies and gives
\[
    K(X,Y) = \frac12\|[X,Y]\|_F^2 = \frac122\|x\times y\|_2^2
    = \|x\|_2^2\|y\|_2^2 \left\|\frac{x}{\|x\|_2}\times \frac{y}{\|y\|_2}\right\|_2^2 =\frac14.
\]
The same reasoning applies to $St(3,3)$.
\end{proof}
Note that with the tools at hand, \eqref{eq:comm_skew3} is a one-liner,
\[
    \|[X,Y]\|_F^2 = 2 \|x\times y\|_2^2 \leq 2 \|x\|_2^2\|y\|_2^2
    =2\frac12\|X\|_F^2\frac12\|Y\|_F^2 =\frac12\|X\|_F^2\|Y\|_F^2.
\]
Equality holds if and only if $X\bot_F Y$.
%
%
%
\section{Lie group essentials}
\label{app:LieGroups}
This section recaps the basics of  Lie groups and quotients of Lie groups. It is mainly collected from the textbooks \cite{gallier2020differential}, \cite{Hall_Lie2015} and  \cite{ONeill1983}.
\subsection{Matrix Lie groups}
\label{app:matLieGroups}
A {\em Lie group} is a differentiable manifold $\mcG$ which
also has a group structure, such that the group operations
`multiplication' and `inversion',
$$\mcG\times
\mcG\ni(g, \tilde{g})\mapsto g\cdot \tilde{g}\in \mcG
\mbox{\hspace{0.1cm} and }\hspace{0.1cm} \mcG\ni g\mapsto g^{-1}\in \mcG
$$ 
are both smooth.
By definition, a {\em matrix Lie group} $\mcG$ is a subgroup of $GL(n,\C)$ that is closed relative to $GL(n,\C)$. It is then also a Lie group in the above sense.

Let $\mcG$ be a real matrix Lie group. When endowed with the bracket operator or {\em matrix commutator} $[X,Y] = XY-YX$, the tangent space $T_I\mcG$ at the identity is the {\em Lie algebra} associated with the Lie group $\mcG$ and is denoted by $\mathfrak{g} = T_I\mcG$.
The linear, skew-symmetric bracket operation is called {\em Lie bracket} and satisfies the {\em Jacobi identity}
\[
[X, [Y,Z]] + [Z,[X,Y]] + [Y,[Z,X]] = 0.
\]
For any $A\in \mcG$, the function ``left-multiplication with $A$'' is a diffeomorphism
	$L_A: \mcG\to \mcG, L_A(B) = AB$; its differential at a point $B\in \mcG$ is the isomporphism
	$$d(L_A)_B: T_B\mcG\to T_{L_A(B)}\mcG, \quad d(L_A)_B(X) = AX.$$
	(Analogous for ``right-multiplication with $A$'', $R_A(B) = BA$, $d(R_A)_B(X) = XA$.)
	The
tangent space at an arbitrary location $A\in \mcG$ is given by the translates (by left-multiplication) of the tangent space at the identity:
\begin{equation}
	\label{eq:tangspaceshifts}
	T_A\mcG = T_{L_A(I)}\mcG =  A \mathfrak{g} = \left\{\Delta = AX\in \R^{n\times n}|\hspace{0.2cm} X\in \mathfrak{g}\right\},
\end{equation}
\cite[\S 5.6, p. 160]{godement2017introduction}.
The Lie algebra $\mathfrak{g} = T_I\mcG$ of $\mcG$ can equivalently be characterized as the set of all matrices $\Delta$ such that $\exp_m(t\Delta)\in \mcG$ for all $t\in \R$, see \cite[Def. 3.18 \& Cor. 3.46]{Hall_Lie2015} for the details.
The exponential map for a matrix Lie group is the matrix exponential restricted to the
corresponding Lie algebra \cite[\S 3.7]{Hall_Lie2015},
\[
\exp_m|_{\mathfrak{g}}: \mathfrak{g}\rightarrow \mcG.
\]
A Riemannian metric $\langle \cdot,\cdot\rangle_g^\mcG$ on $\mcG$ is called
 {\em left-invariant} if 
 $$\langle d(L_A)_B X, d(L_A)_B Y \rangle_{AB}^\mcG (=\langle A X, A Y \rangle_{AB}^\mcG) = \langle X,Y\rangle_B^\mcG $$
 for all $X,Y\in T_B\mcG$.
 It is called right-invariant, if $$\langle d(R_A)_B X, d(R_A)_B Y \rangle_{BA}^\mcG (=\langle XA,  YA\rangle_{BA}^\mcG )= \langle X,Y\rangle_B^\mcG $$ for all $X,Y\in T_B\mcG$.
 If a metric is left- and right-invariant, it is called {\em bi-invariant}.
\subsection{Quotients of Lie groups by closed subgroups}
		Let $\mcG$ be a Lie group and $\mcH\leq \mcG$ be a Lie subgroup. 
		For $A\in \mcG$, a subset of $\mcG$ of the form 
		$
		[A]:= A\mcH = \{A\cdot Q|\hspace{0.1cm} Q\in \mcH\}
		$
		is called a {\em left coset of $\mcH$}.
		The left coset $[I]$ is the subgroup itself.
		The left cosets form a partition of $\mcG$,
		and the quotient space $\mcG/\mcH$ determined by this partition is called the {\em left coset space of $\mcG$ modulo $\mcH$}, see \cite[\S 21, p. 551]{Lee2012smooth}.
	The next is the central theorem for quotients of Lie groups.
	\begin{theorem}
		\label{thm:HomSpaceConst}
		(cf. \cite[Thm 21.17, p. 551]{Lee2012smooth})
		Let $\mcG$ be a Lie group and let $\mcH$ be a closed subgroup of $\mcG$. The left coset space $\mcG/\mcH$ is a manifold of dimension $\dim \mcG - \dim \mcH$ with a unique differentiable structure such that the quotient map 
		$\pi: \mcG \to  \mcG/\mcH, A\mapsto [A]$ is a smooth surjective submersion.
		The left action of $\mcG$ on $\mcG/\mcH$ given by $A(B \mcH) = (AB) \mcH$
		turns $\mcG/\mcH$ into a {\em homogeneous $\mcG$-space.}
	\end{theorem}

    Each preimage $\mcG_{\pi(A)}:=\pi^{-1}([A]) \subset \mcG$, called fiber, is a closed, embedded submanifold.
	Under the Riemannian metric $\langle \cdot,\cdot\rangle_A^\mcG$, 
	at each point $A\in \mcG$, the tangent space $T_A\mcG$
	decomposes into an {\em orthogonal  direct sum} $T_A\mcG = T_A\mcG_{\pi(A)}\oplus (T_A\mcG_{\pi(A)})^\bot$
	with respect to the metric.
	The tangent space of the fiber $T_A\mcG_{\pi(A)}=: V_A$ 
	 is the kernel of the differential 
	$d\pi_A: T_A\mcG \to T_{\pi(A)}\mcG/\mcH$ and is called the {\em vertical space}.
	Its orthogonal complement $H_A := V_A^{\bot}$ is the {\em horizontal space}.
	The key issue is that the tangent space of the quotient at $[A] = \pi(A)$ may be identified with the horizontal space at $A$.
	 
	\[
		H_A \cong T_{[A]}\mcG/\mcH.
    \]
\begin{itemize} 
	\item 
	For every tangent vector $Y\in T_{[A]}(\mcG/\mcH)$ there is $\bar Z =\bar X+\bar Y\in V_A \oplus H_A = T_A\mcG$ such that 	$d\pi_A(\bar Z) = Y$. The horizontal component $\bar Y\in H_A$ is unique and is called the horizontal lift of $Y\in T_{\pi(A)}(\mcG/\mcH)$.
	By going to the horizontal lifts, a Riemannian metric on the quotient can be defined by
	\begin{equation}
		\label{eq:quotMetric}
		\langle Y_1,Y_2 \rangle_{[A]}^{\mcG/\mcH}
		:= \langle \bar Y_1, \bar Y_2 \rangle_A^{\mcG} 
	\end{equation}
	for $Y_1,Y_2 \in T_{[A])}(\mcG/\mcH)$. 
	\item   W.r.t. this (and only this) metric, by construction,
	$d\pi_A$ preserves inner products:
	$d\pi_A$ is a linear isometry
	between $H_A$ and  $T_{[A]}(\mcG/\mcH)$.
	
	\item Horizontal geodesics in $\mcG$ are mapped to geodesics in $\mcG/\mcH$ under $\pi$.
	Horizontal geodesics are geodesics with velocity field staying in the horizontal space for all time $t$.
\end{itemize}
At the special point $A=I$, the vertical space is the Lie algebra of $\mcH$,  $V_I = \ker{d\pi_I} = \mathfrak{h}$.
This is because $\mcH=\pi^{-1}(I)$. For any curve $C(t)\subset \pi^{-1}(I)=\mcH$ starting from $C(0)=I$ 
with image in the fiber, it holds $\dot C (0)\in T_I\mcH=\mathfrak h$. On the other  hand, $\pi$ is constant along the fiber so that $0=\frac{d}{dt}\big\vert_{t=0} \pi(C(t)) = d\pi_I[\dot C(0)]$.
Hence, any vector tangent to $\mcH$ at $I$ is in the kernel of $d\pi_I$. 
At the identity $I$, the splitting into vertical and horizontal space is
\[
	T_I\mathcal{G} = \mathfrak h \oplus \mathfrak m = V_I\oplus H_I = T_I\mcH \oplus (T_I\mcH)^\bot.
\]
Hence, the tangent space of the quotient at $\pi(I)$ is $\mathfrak m \cong T_{\pi(I)}\mcG/\mcH$.
This choice of symbols is common in the literature, but the reader should be aware that $\mathfrak h$ is the vertical space, with the choice of letter referring to the subgroup name $\mcH$ and not to "h for horizontal".
The choice of the symbol $\mathfrak m$ is motivated by the fact that if the quotient space is called
$\mcM:=\mcG/\mcH$, then $\mathfrak m$ is the tangent space at $\pi(I)$.
It is not the associated Lie algebra though, because in general $\mcG/\mcH$ is not a Lie group.

\end{document}